\documentclass[A4]{article}
\usepackage[latin1]{inputenc}
\usepackage[english]{babel}
\usepackage{amsmath,amssymb,amsthm}
\usepackage{psfrag}
\usepackage[]{mathrsfs}
\usepackage{graphicx}
\usepackage[all]{xypic}

\newcommand{\asv}[1]{v_{#1}}
\newcommand{\upartial}{\mathrm{\partial}}

\newtheorem{lemma}{Lemma}
\newtheorem{definition}{Definition}
\newtheorem{theorem}{Theorem}
\newtheorem{proposition}{Proposition}
\newtheorem{corollary}{Corollary}
\newtheorem{remark}{Remark}
\title{An entire transcendental family with a persistent Siegel disc}
\author{Rubén Berenguel\thanks{Email: ruben@maia.ub.es}, Núria Fagella\thanks{Email: fagella@maia.ub.es}}

\usepackage{enumerate}
\DeclareMathAlphabet{\curv}{OT1}{pzc}{m}{it}
\begin{document}

\maketitle
\begin{abstract}
  We study the class of entire transcendental maps of finite order
  with one critical point and one asymptotic value, which has exactly
  one finite pre-image, and having a persistent Siegel disc. After
  normalisation this is a one parameter family
  $f_a$ with $a\in\mathbb{C}^*$ which includes the semi-standard map
  $\lambda ze^z$ at $a=1$, approaches the exponential map when $a\to0$
  and a quadratic polynomial when $a\to\infty$. We investigate the
  stable components of the parameter plane (capture components and
  semi-hyperbolic components) and also some topological properties of
  the Siegel disc in terms of the parameter.
\end{abstract}

\section{Introduction}
Given a holomorphic endomorphism $f:S\to S$ on a Riemann surface $S$
we consider the dynamical system generated by the iterates of $f$,
denoted by $f^n=f\circ\stackrel{n}{\cdots}\circ f$. The {
\em orbit} of an initial condition $z_0 \in S$ is the sequence
$\mathcal{O}^+(z_0)=\{f^n(z_0)\}_{n\in \mathbb{N}}$ and we are interested in
classifying the initial conditions in the {\em phase space} or {\em
dynamical plane} $S$, according to the asymptotic behaviour of their
orbits when $n$ tends to infinity.

There is a dynamically natural partition of the phase space $S$ into
the {\em Fatou set} $\curv{F}(f)$ (open) where the iterates of $f$
form a normal family and the {\em Julia set} $\curv{J}(f)=S
\backslash\curv{F}(f)$ which is its complement (closed).

If $S=\widehat{\mathbb{C}}=\mathbb{C} \cup \infty$ then $f$ is a
rational map. If $S=\mathbb{C}$ and $f$ does not extend to the point at
infinity, then $f$ is an {\em entire transcendental map}, that is,
infinity is an essential singularity. Entire transcendental functions
present many differences with respect to rational maps.

One of them concerns the singularities of the inverse function. For a
rational map, all branches of the inverse function are well defined
except at a finite number of points called the {\em critical values},
points $w=f(c)$ where $f'(c)=0$. The point $c$ is then called a {\em
critical point}. If $f$ is an entire transcendental map, there is
another possible obstruction for a branch of the inverse to be well
defined, namely its {\em asymptotic values}. A point $v\in \mathbb{C}$ is
called an asymptotic value if there exists a path $\gamma(t) \to
\infty$ when $t\to\infty$, such that $f(\gamma(t))\to v$ as $t\to
\infty$.  An example is $v=0$ for $f(z)=e^z$, where $\gamma(t)$ can be
chosen to be the negative real axis.

In any case, the set of singularities of the inverse function, also
called {\em singular values}, plays a very important role in the
theory of iteration of holomorphic functions. This statement is
motivated by the non-trivial fact that most connected components of the
Fatou set (or stable set) are somehow associated to a singular
value. Therefore, knowing the behaviour of the singular orbits provides
information about the nature of the stable orbits in the phase space.

The dynamics of rational maps are fairly well understood, given the
fact that they possess a finite number of critical points and hence of
singular values. This motivated the definition and study of special
classes of entire transcendental functions like, for example, the
class $\mathcal{S}$ of functions of {\em finite type} which are those
with a finite number of singular values. A larger class is
$\mathcal{B}$ the class of functions with a bounded set of
singularities. These functions share many properties with rational
maps, one of the most important is the fact that every connected
component of the Fatou set is eventually periodic (see e.g.
\cite{EremenkoLyubichDynamicalEntire92} or
\cite{GoldbergKeenFinitenessEntire86}). There is a classification of
all possible periodic connected components of the Fatou set for
rational maps or for entire transcendental maps in class
$\mathcal{S}$. Such a component can only be part of a cycle of
rotation domains (Siegel discs) or part of the basin of attraction of
an attracting, super-attracting or parabolic periodic orbit.

We are specially interested in the case of rotation domains. We say
that $\Delta$ is an invariant Siegel disc if there exists a conformal
isomorphism $\varphi:\Delta \to \mathbb{D}$ which conjugates $f$ to
$\mathcal{R}_\theta(z)=e^{2\pi i\theta} z$ (and $\varphi$ can not be
extended further), with $\theta \in \mathbb{R}\setminus
\mathbb{Q}\cap(0,1)$ called the {\em rotation number} of $\Delta$.
Therefore a Siegel disc is foliated by invariant closed simple curves,
where orbits are dense. The existence of such Fatou components was
first settled by Siegel \cite{SiegelIteration} who showed that if
$z_0$ is a fixed point of multiplier $\rho=f'(z_0)=e^{2\pi i \theta}$
and $\theta$ satisfies a Diophantine condition, then $z_0$ is {\em
  analytically linearisable} in a neighbourhood or, equivalently, $z_0$
is the centre of a Siegel disc. The Diophantine condition was relaxed
later by Brjuno and R\"ussman (for an account of these proofs see e.g.
\cite{MilnorDynamics}), who showed that the same is true if $\theta$
belonged to the set of Brjuno numbers $\mathcal{B}$. The relation of
Siegel discs with singular orbits is as follows. Clearly $\Delta$
cannot contain critical points since the map is univalent in the disc.
Instead, the boundary of $\Delta$ must be contained in the {\em
  post-critical set }
$\cup_{c\in\mathrm{Sing}(f^{-1})}\overline{\mathcal{O}^+(c)}$ i.e.,
the accumulation set of all singular orbits. In fact something
stronger is true, namely that $\partial \Delta$ is contained in the
accumulation set of the orbit of at least {\em one} singular value
(see \cite{ManeFatouSiegel}).

Our goal in this paper is to describe the dynamics of the one
parameter family of entire transcendental maps
\[ 
f_a(z) = \lambda a(e^{z/a}(z+1-a)-1+a),
\] where $a\in \mathbb{C}\setminus \{0\}=\mathbb{C}^*$ and $\lambda=e^{2\pi
i\theta}$ with $\theta$ being a fixed irrational Brjuno number.
Observe that $0$ is a fixed point of multiplier $\lambda$ and
therefore, for all values of the parameter $a$, there is a persistent
Siegel disc $\Delta_a$ around $z=0$. The functions $f_a$ have two
singular values: the image of the only critical point $w=-1$ and an
asymptotic value at $\asv{a}=\lambda a (a-1)$ which has one and only
one finite pre-image at the point $p_a=a-1$.

The motivation for studying this family of maps is manifold. On one
hand this is the simplest family of entire transcendental maps having
one simple critical point and one asymptotic value with a finite
pre-image (see Theorem \ref{TheoremCharacterizationOfTheFamily} for
the actual characterisation of $f_a$). The persistent Siegel disc
makes it into a one-parameter family, since one of the two singular
orbits must be accumulating on the boundary of $\Delta_a$. We will see
that the situation is very different, depending on which of the two
singular values is doing that. Therefore, these maps could be viewed
as the transcendental version of cubic polynomials with a persistent
invariant Siegel disc, studied by Zakeri in \cite{ZakeriCubicSiegel}.
In our case, many new phenomena are possible with respect to the cubic
situation, like unbounded Siegel discs for example; but still the two
parameter planes share many features like the existence of capture
components or semi-hyperbolic ones.

There is a second motivation for studying the maps $f_a$, namely that
this one parameter family includes in some sense three emblematic
examples. For $a=1$ we have the function $f_1(z)=\lambda z e^z$, for
large values of $a$ we will see that $f_a$ is polynomial-like of
degree 2 in a neighbourhood of the origin (see Theorem
\ref{TheoremTheFamilyIsPolynomialLike}); finally when $a\to 0$, the
dynamics of $f_a$ are approaching those of the exponential map
$u\mapsto \lambda(e^u-1)$, as it can be seen changing variables to
$u=z/a$. Thus the parameter plane of $f_a$ can be thought of as
containing the polynomial $\lambda(z +\frac{z^2}{2})$ at infinity, its
transcendental analogue $f_1$ at $a=1$, and the exponential map at
$a=0$. The maps $z\mapsto \lambda z e^z$ have been widely studied (see
\cite{GeyerSiegelHerman} and \cite{FagellaLimitingDynamics}), among
other reasons, because they share many properties with quadratic
polynomials: in particular it is known that when $\theta$ is of
constant type, the boundary of the Siegel disc is a quasi-circle that
contains the critical point. It is not known however whether there
exist values of $\theta$ for which the Siegel disc of $f_1$ is
unbounded. In the long term we hope that this family $f_a$ can throw
some light into this and other problems about $f_1$.

For the maps at hand we prove the following.
\renewcommand{\thetheorem}{A}
\begin{theorem}
  \begin{enumerate}[a)]
  \item \label{PartATmaA}There exists $R,M>0$ such that if $\theta$ is
    of constant type and $|a|>M$ then the boundary of $\Delta_a$ is a
    quasi-circle which contains the critical point. Moreover
    $\Delta_a\subset D(0,R)$.
  \item \label{PartBTmaA} If $\theta$ is Diophantine and the orbit of
    $c=-1$ belongs to a periodic basin or is eventually captured by
    the Siegel disc, then either the Siegel disc $\Delta_a$ is
    unbounded or its boundary is an indecomposable continuum.
  \item \label{PartCTmaA} If $\theta$ is Diophantine and
    $f_a^n(-1)\stackrel{n\to\infty}{\longrightarrow}\infty$ the Siegel
    disc $\Delta_a$ is unbounded, and the boundary contains the
    asymptotic value.
  \end{enumerate}
\end{theorem}

Part \ref{PartATmaA}) follows from Theorem
\ref{TheoremTheFamilyIsPolynomialLike} (see Corollary
\ref{CorollaryConstantTypeQuasiCircle} below it). The remaining parts
(Theorem \ref{ThmUnboundedSiegelDisks}) are based on Herman's proof
\cite{HermanCriticalExponential} of the fact that Siegel discs of the
exponential map are unbounded, if the rotation number is Diophantine,
although in this case there are some extra difficulties given by the
free critical point and the finite pre-image of the asymptotic value.

In this paper we are also interested in studying the parameter plane of
$f_a$, which is $\mathbb{C}^*$, and in particular the connected
components of its {\em stable} set, i.e., the parameter values for
which the iterates of both singular values form a normal family in
some neighbourhood. We denote this set as $\mathcal S$ (not to be
confused with the class of finite type functions). These connected
components are either {\em capture components}, where an iterate of
the free singular value falls into the Siegel disc; or {\em
semi-hyperbolic}, when there exists an attracting periodic orbit
(which must then attract the free singular value); otherwise they are
called {\em queer}.


The following theorem summarises the properties of semi-hyperbolic
components, and is proved in Section
\ref{SectionSemiHyperbolicComponents} (see Proposition
\ref{PropositionHyperbolicSimplyConnected}, Theorems
\ref{TheoremH1vIsUnbounded},
\ref{TheoremParametrizationSemiHyperbolicAsymptotic} and Proposition
\ref{PropositionParametrizationSemiHyperbolicCritical} therein). By a
\emph{component} of a set we mean a connected component.
\renewcommand{\thetheorem}{B}
\begin{theorem}
  Define 
  \begin{align}
    H^c&=\{a\in\mathbb{C}\vert\mathcal{O}^+(-1)\textrm{ is attracted
      to
      an attracting periodic orbit}\},\nonumber\\
    H^v&=\{a\in\mathbb{C}\vert\mathcal{O}^+(\asv{a})\textrm{ is
      attracted to an attracting periodic orbit}\}.\nonumber
  \end{align}
\begin{enumerate}[a)]
\item \label{PartATmaB} Every component of $H^v\cup H^c$ is simply
  connected.
\item \label{PartBTmaB} If $W$ is a component of $H^v$ then $W$ is
  unbounded and the multiplier map $\chi:W\to\mathbb{D}^*$ is the
  universal covering map.
\item \label{PartCTmaB} There is one component $H_1^v$ of $H^v$ for
  which $\mathcal{O}^+(\asv{a})$ tends to an attracting fixed point.
  $H_1^v$ contains the segment $[r,\infty)$ for $r$ large enough.
\item \label{PartDTmaB} If $W$ is a component of $H^c$, then $W$ is
  bounded and the multiplier map $\chi:W\to\mathbb{D}$ is a conformal
  isomorphism.
\end{enumerate}
\end{theorem}

Indeed, when the critical point is attracted by a cycle, we naturally
see copies of the Mandelbrot set in parameter space. Instead, when it
is the asymptotic value that acts in a hyperbolic fashion, we find
unbounded exponential-like components, which can be parametrised using
quasi-conformal surgery.

\begin{figure}[!hbt]
  \begin{center}
    \psfrag{2}{\small$H_2$} \psfrag{1}{\small$H_1$}
    \psfrag{0}{\small$0$} \psfrag{C}{\small$1$}
    \fbox{\includegraphics[width=0.45\textwidth]{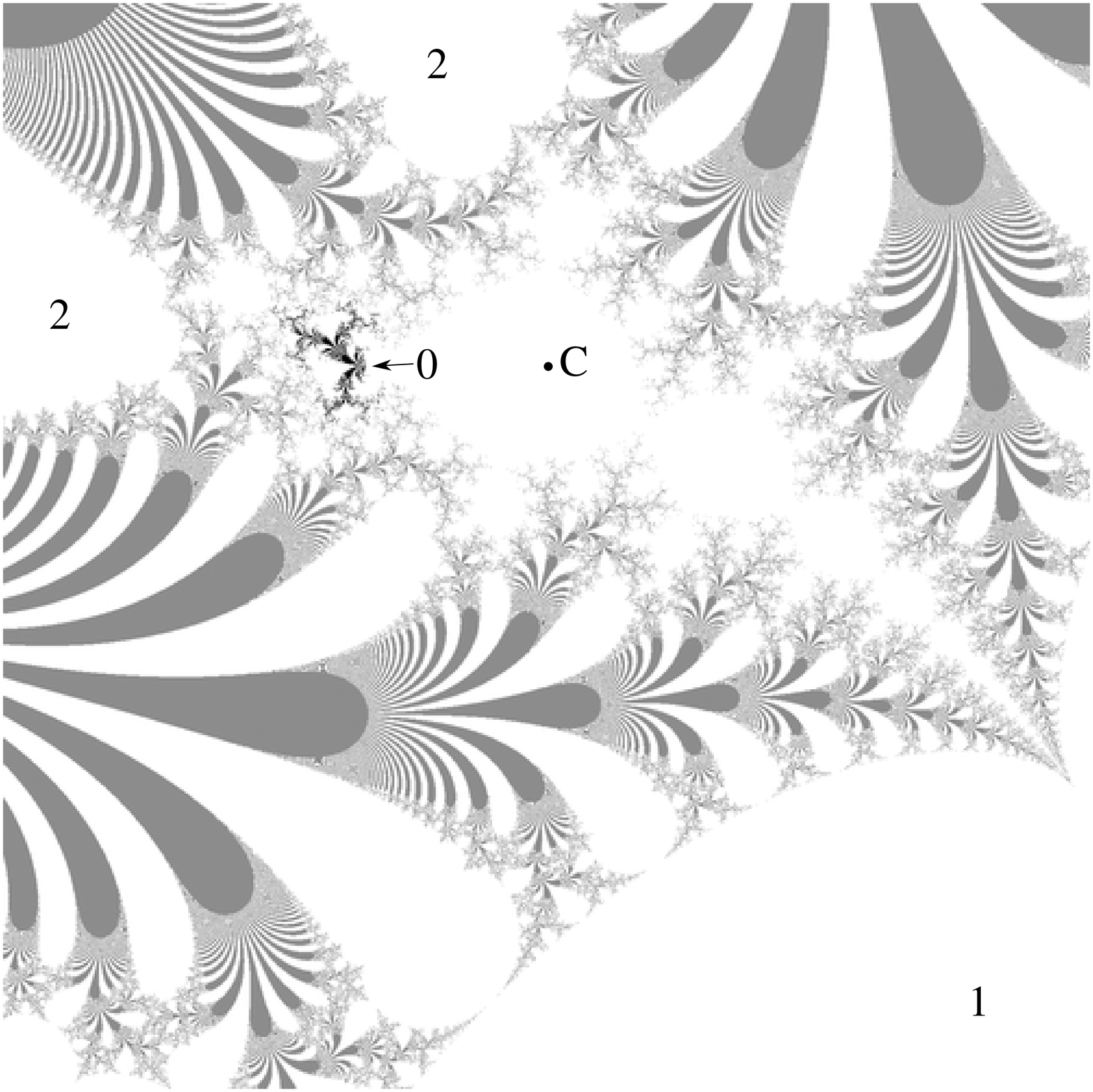}
    }
    \fbox{\includegraphics[width=0.45\textwidth]{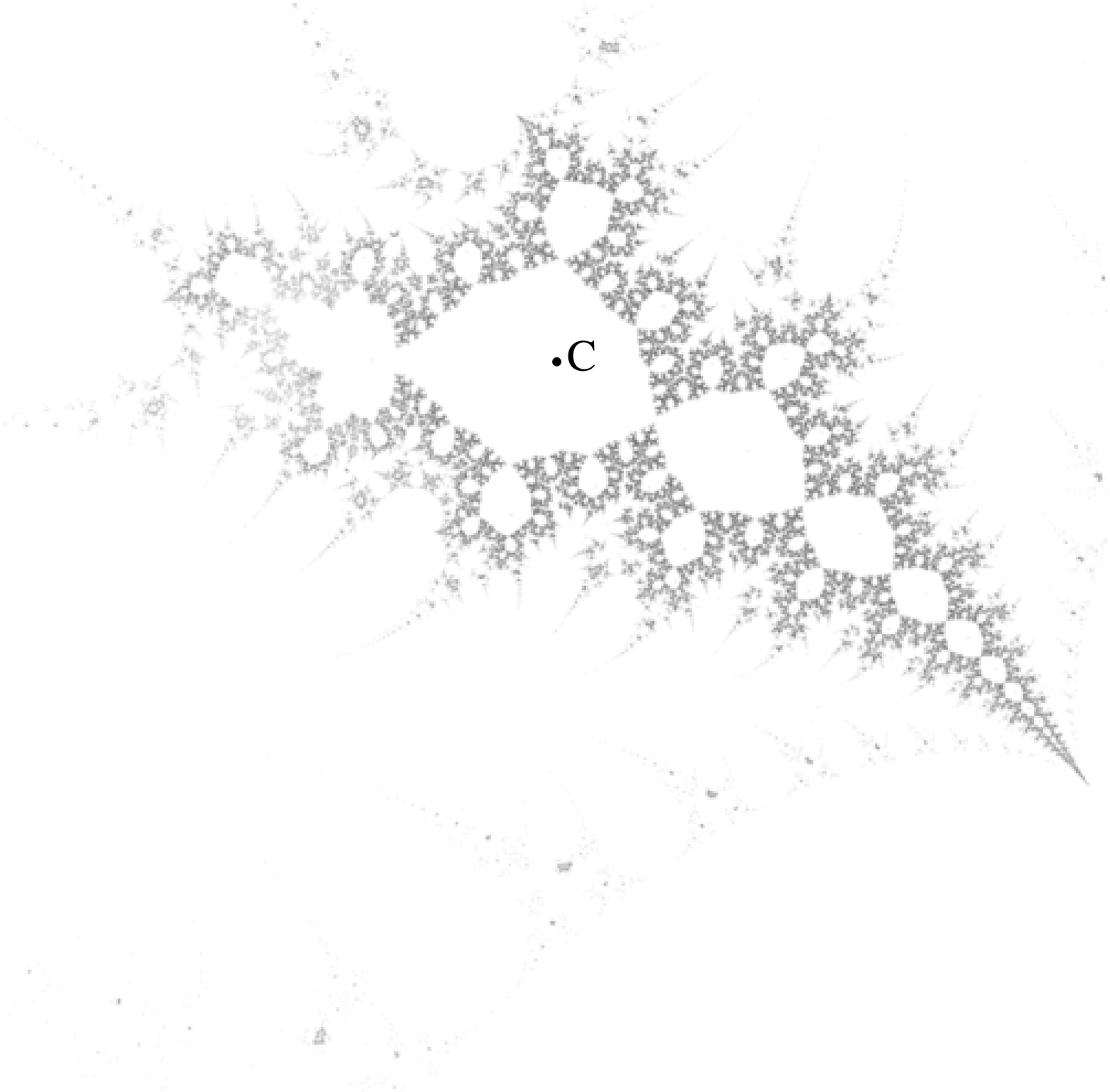}
    }
    \caption{{\bf Left:} Simple escape time plot of the parameter
      plane. Light grey: asymptotic orbit escapes, dark grey critical
      orbit escapes, white neither escapes. Regions labelled $H_1$ and
      $H_2$ correspond to parameters for which the asymptotic value is
      attracted to an attracting cycle. {\bf Right:} The same plot,
      using a different algorithm which emphasises the capture
      components. Upper left: $(-2,2),$ Lower right: $(4,-4)$.}
  \label{FigureParametersEscapeToInfinityZoomOut}
  \end{center}
\end{figure}

A dichotomy also occurs with capture components. Numerically we can
observe copies of quadratic Siegel discs in parameter space, which
correspond to components for which the asymptotic value is being
captured. There is in fact a main capture component $C_0^v$, the one
containing $a=1$ (see Figure
\ref{FigureParametersEscapeToInfinityZoomOut}), which corresponds to
parameters for which the asymptotic value $v_a$, belongs itself to the
Siegel disc. This is possible because of the existence of a finite
pre-image of $v_a$. The centre of $C_0^v$ is the semi-standard map
$f_1(z)=\lambda ze^z$, for which zero itself is the asymptotic value.

The properties we show for capture components are summarised in the
following theorem (see Section \ref{SectionCaptureComponents}: Theorem
\ref{TheoremCaptureOpenCondition} and Proposition
\ref{PropositionCaptureSimplyConnected}).

\renewcommand{\thetheorem}{C}
\begin{theorem}
  Let us define
  \begin{align}
    C^c&=\{a\in\mathbb{C}\vert f_a^n(-1)\in\Delta_a\textrm{ for some }n\geq1\},\nonumber\\
    C^v&=\{a\in\mathbb{C}\vert f_a^n(\asv{a})\in\Delta_a\textrm{ for
      some }n\geq0\}.\nonumber
   \end{align}
   Then
\begin{enumerate}[a)]
\item \label{PartATmaC} $C^c$ and $C^v$ are open sets.
\item \label{PartBTmaC} Every component $W$ of $C^c\cup C^v$ is simply
  connected.
\item \label{PartCTmaC} Every component $W$ of $C^c$ is bounded.
\item \label{PartDTmaC} There is only one component of $C_0^v=\{a\in
  \mathbb{C}\vert\asv{a}\in\Delta_a\}$ and it is bounded.
\end{enumerate}
\end{theorem}

Numerical experiments show that if $\theta$ is of constant type, the
boundary of $C_0^v$ is a Jordan curve, corresponding to those parameter
values for which both singular values lie on the boundary of the
Siegel disc (see Figure
\ref{FigureParametersEscapeToInfinityZoomOut}). This is true for the
slice of cubic polynomials having a Siegel disc of rotation number
$\theta$, as shown by Zakeri in \cite{ZakeriCubicSiegel}, but his
techniques do not apply to this transcendental case.

As we already mentioned, we are also interested in parameter values
for which $f_a$ is Julia stable, i.e. where both families of iterates
$\{f_a^n(-1)\}_{n\in\mathbb{N}}$ and $\{f_a^n(\asv{a})\}_{n\in\mathbb{N}}$
are normal in a neighbourhood of $a$ (see Section
\ref{SectionJuliaStability}). We first show that any parameter in a
capture component or a semi-hyperbolic component is
$\mathcal{J}$-stable.
\renewcommand{\thetheorem}{D}
\begin{proposition}
  If $a\in H\cup C$ then $f_a$ is $\mathcal{J}$-stable, where
  $H=H^c\cup H^v$ and $C=C^c\cup C^v$.
\end{proposition}

By using holomorphic motions and the proposition above, it is enough
to have certain properties for one parameter value $a_0$, to be able
to ``extend'' them to all parameters belonging to the same stable
component. More precisely we obtain the following corollaries (see
Proposition \ref{PropositionHolomorphicMotionOfC0} and Corollary
\ref{CorollaryPolynomialLikeConsequencesHolomorphicMotions}).
\renewcommand{\thetheorem}{E}
\begin{proposition}
   \begin{enumerate}[a)]
   \item \label{PartAPropositionE} If $\theta$ is of constant type and
     $a\in C_0^v$ (i.e. the asymptotic value lies inside the Siegel
     disc) then $\partial\Delta_a$ is a quasi-circle that contains the
     critical point.
   \item \label{PartBPropositionE} Let $W\subset H^v\cup C^v$ be a
     component intersecting $\{|z|>M\}$ where $M$ is as in Theorem
     A. Then,
     \begin{enumerate}[i)]
     \item if $\theta$ is of constant type, for all $a\in W$ the
       boundary $\partial\Delta_a$ is a quasi-circle containing the
       critical point.
     \item There exist values of
       $\theta\in\mathbb{R}\backslash\mathbb{Q}\cap(0,1)$ such that if
       a component $W\subset C^v\cup H^v$ intersects $\{|z|>M\}$, then
       for all $a\in W$, the boundary of $\Delta_a$ is a quasi-circle
       \emph{not} containing the critical point.  
     \end{enumerate}
   \end{enumerate}
\end{proposition}




The paper is organised as follows. Section \ref{SectionPreliminars}
contains statements and references of some of the results used
throughout the paper. Section
\ref{SectionCharacterizationImagesDynamicalPlanes} contains the
characterisation of the family $f_a$, together with descriptions and
images of the possible scenarios in dynamical plane. It also
contains the proof of Theorem A. Section
\ref{SectionSemiHyperbolicComponents} deals with semi-hyperbolic
components and contains the proof of Theorem B, split in several
parts, and not necessarily in order. In the same fashion, capture
components and Theorem C are treated in Section
\ref{SectionCaptureComponents}. Finally Section
\ref{SectionJuliaStability} investigates Julia stability and contains
the proofs of Propositions D and E.
 
\renewcommand\thetheorem{\thesection.\arabic{theorem}}
 

\section{Preliminary results}
\label{SectionPreliminars}
In this section we state results and definitions which will be useful
in the sections to follow.

\subsection{Quasi-conformal mappings and holomorphic motions}

First we introduce the concept of \emph{quasi-conformal mapping}.
Quasi-conformal mappings are a very useful tool in complex dynamical
systems as they provide a bridge between a geometric construction for
a system and its analytic information. They are also a fundamental
pillar for the framework of \emph{quasi-conformal surgery}, the other
one being the \emph{measurable Riemann mapping theorem}. For the
groundwork on quasi-conformal mappings see for example
\cite{AhlforsLQC2006}, and for an exhaustive account on
quasi-conformal surgery, see \cite{BrannerFagellaQCS}.

\begin{definition}
  Let $\mu:U\subseteq\mathbb{C}\to\mathbb{C}$ be a measurable
  function. Then it is a $k$-Beltrami form (or Beltrami coefficient,
  or complex dilatation) of $U$ if $\|\mu(z)\|_\infty\leq k<1$.
\end{definition}

\begin{definition}
  Let $f:U\subseteq\mathbb{C}\to V\subseteq\mathbb{C}$ be a
  homeomorphism. We call it $k$-quasi-conformal if locally it has
  distributional derivatives in $\mathcal{L}^2$ and
  \begin{align}
    \mu_f(z)=\frac{\frac{\upartial f}{\upartial\bar
        z}(z)}{\frac{\upartial f}{\upartial z}(z)}
  \end{align}
  is a $k$-Beltrami coefficient. Then $\mu_f$ is called the complex
  dilatation of $f(z)$ (or the Beltrami coefficient of $f(z)$).

  Given $f(z)$ satisfying all above except being an homeomorphism, we
  call it $k$-quasi-regular.
\end{definition}

The following technical theorem will be used when we have compositions
of quasi-conformal mappings and finite order mappings.

\begin{theorem}[{\cite[p.~750]{FagellaSearaHerman04}}]\label{TheoremAhlforsHolderContinuity}
  A $k$-quasi-conformal mapping in a domain $U\subset\mathbb{C}$ is
  uniformly Hölder continuous with exponent $(1-k)/(1+k)$ in every
  compact subset of $U$.
\end{theorem}

\begin{theorem}[(Measurable Riemann Mapping, MRMT)]
  Let $\mu$ be a Beltrami form over $\mathbb{C}$.  Then there exists a
  quasi-conformal homeomorphism $f$ integrating $\mu$ (i.e. the
  Beltrami coefficient of $f$ is $\mu$), unique up to composition with
  an affine transformation.
\end{theorem}
\begin{theorem}[(MRMT with dependence of parameters)]
  Let $\Lambda$ be an open set of $\mathbb{C}$ and let
  $\{\mu_\lambda\}_{\lambda\in\Lambda}$ be a family of Beltrami forms
  on $\mathbb{\hat C}$.  Suppose $\lambda\to\mu_\lambda(z)$ is
  holomorphic for each fixed $z\in\mathbb{C}$ and
  $\|\mu_\lambda\|_\infty\leq k<1$ for all $\lambda$. Let $f_\lambda$
  be the unique quasi-conformal homeomorphism which integrates
  $\mu_\lambda$ and fixes three given points in $\mathbb{\hat
    C}$. Then for each $z\in\mathbb{\hat C}$ the map $\lambda\to
  f_{\lambda}(z)$ is holomorphic.
\end{theorem}

The concept of holomorphic motion was introduced in
\cite{ManeSadSullivanHolMot} along with the (first) $\lambda$-lemma.

\begin{definition}
  Let $h:\Lambda\times X_0\to \mathbb{\hat C},$ where $\Lambda$
  is a complex manifold and $X_0$ an arbitrary subset of
  $\mathbb{\hat C}$, such that
\begin{itemize}
\item $h(0,z)=z$,
\item $h(\lambda,\cdot)$ is an injection from $X_0$ to $\mathbb{\hat C}$,
\item For all $z\in X_0$, $z\mapsto h(\lambda,z)$ is holomorphic. 
\end{itemize}
Then $h_\lambda(z)=h(\lambda,z)$ is called a \emph{holomorphic motion}
of $X$.
\end{definition}

The following two fundamental results can be found in
\cite{ManeSadSullivanHolMot} and \cite{SlodkowskiHolMot} respectively.

\begin{lemma}[(First $\lambda$-lemma)]
  A holomorphic motion $h_\lambda$ of any set
  $X\subset\mathbb{\hat C}$ extends to a jointly continuous
  holomorphic motion of $\bar X$.
\end{lemma}

\begin{lemma}[(Second $\lambda$-lemma)]\label{LemmaSecondLambda}
  Let $U\subset \mathbb{C}$ be a set and $h_\lambda$ a holomorphic
  motion of $U$. This motion extends to a holomorphic motion of
  $\mathbb{C}$.
\end{lemma}

\subsection{Hadamard's factorisation theorem}

We will need the notion of \emph{rank} and \emph{order} to be able to
state Hadamard's factorisation theorem, which we will use in the proof
of Theorem \ref{TheoremCharacterizationOfTheFamily}. All these results
can be found in \cite{ConwayFunctions1}.

\begin{definition}
  Given $f:\mathbb{C}\to\mathbb{C}$ an entire function we say it is of
  \emph{finite order} if there are positive constants $a>0$, $r_0>0$
  such that
  \[|f(z)|<e^{|z|^a},\quad \text{for }|z|>r_0.\] Otherwise, we say
  $f(z)$ is of infinite order.  We define
  \[\lambda=\inf\{a\vert |f(z)|<\mathrm{exp}(|z|^a)\text{ for }
  |z|\text{ large enough}\}\] as the order of $f(z)$.
\end{definition}

\begin{definition}
  Let $f:\mathbb{C}\to\mathbb{C}$ be an entire function with zeroes
  $\{a_1,a_2,\ldots\}$ counted according to multiplicity. We say $f$
  is of \emph{finite rank} if there is an integer $p$ such that
  \begin{align}\sum_{n=1}^\infty|a_n|^{p+1}<\infty.\label{EquationFiniteRank}\end{align}
  We say it is of rank $p$ if $p$ is the smallest integer verifying
  \eqref{EquationFiniteRank}. If $f$ has a finite number of zeroes then it has
  rank 0 by definition.
\end{definition}

\begin{definition}
  An entire function $f:\mathbb{C}\to\mathbb{C}$ is said to be of
  \emph{finite genus} if it has finite rank $p$ and it factorises as:
  \begin{align}
    f(z)=z^me^{g(z)}\cdot\prod_{n=1}^\infty E_p(z/a_n),
    \label{EquationHadamardFactorization}
  \end{align}
  where $g(z)$ is a polynomial, $a_n$ are the zeroes of $f(z)$ as in
  the previous definition and
  \[E_p(z)=(1-z)e^{z+\frac{z^2}{2}+\ldots+\frac{z^p}{p}}.\] We define
  the genus of $f(z)$ as $\mu=\max\{\mathrm{deg}\,g, \mathrm{rank}\,f\}$
\end{definition}

\begin{theorem}
  If $f$ is an entire function of finite genus $\mu$ then $f$ is of
  finite order $\lambda<\mu+1$.
\end{theorem}

The converse of this theorem is also true, as we see below.

\begin{theorem}[(Hadamard's factorisation)]\label{TheoremHadamardsFactorization}
  Let $f$ be an entire function of finite order $\lambda$. Then $f$ is
  of finite genus $\mu\leq\lambda$.
\end{theorem}

Observe that Hadamard's factorisation theorem implies that every
entire function of finite order can be factorised as in
\eqref{EquationHadamardFactorization}.





\subsection{Siegel discs}

The following theorem (which is an extension of the original theorem
by C.L. Siegel) gives arithmetic conditions on the rotation number of a
fixed point to ensure the existence of a Siegel disc around it. J-C.
Yoccoz proved that this condition is sharp in the quadratic family.
The proof of this theorem can be found in \cite{MilnorDynamics}.

\begin{theorem}[(Brjuno-Rüssmann)]\label{TheoremLinearizationOfBrjunoRussmann}
Let $f(z)=\lambda z+\mathcal{O}(z^2)$. 
If $\frac{p_n}{q_n}=[a_1;a_2,\ldots,a_n]$ is the $n$-th convergent of
the continued fraction expansion of
$\theta$, where $\lambda=e^{2\pi i\theta}$, and
\begin{align}
\label{EquationBrjunoCondition}
\sum_{n=0}^\infty\frac{\log(q_{n+1})}{q_n}<\infty,\end{align}
then $f$
is locally linearisable.
\end{theorem}

Irrational numbers with this property are called of \emph{Brjuno type}.

We define the notion of \emph{conformal capacity} as a measure of the
``size'' of Siegel discs.

\begin{definition}
  Consider the Siegel disc $\Delta$ and the unique linearising map
  $h:\mathbb{D}(0,r)\stackrel{\sim}{\to}\Delta$, with $h(0)$ and
  $h'(0)=1$. The radius $r>0$ of the domain of $h$ is called the
  \emph{conformal capacity} of $\Delta$ and is denoted by
  $\kappa(\Delta)$.
\end{definition}

A Siegel disc of capacity $r$ contains a disc of radius $\frac{r}{4}$
by Koebe 1/4 Theorem.

The following theorem (see \cite{YoccozDiviseurs} for a proof) shows
that Siegel discs can not shrink indefinitely.

\begin{theorem}\label{TheoremSiegelDisksDoNotShrinkYoccoz}
  Let $0<\theta<1$ be an irrational number of Brjuno type, and let
  $\Phi(\theta)=\sum_{n=1}^\infty(\log q_{n+1}/q_n)<\infty$ be the
  Brjuno function. Let $S(\theta)$ be the space of all univalent
  functions $f:\mathbb{D}\to\mathbb{C}$ with $f(0)=0$ and
  $f'(0)=e^{2\pi i\theta}$. Finally, define $\kappa(\theta)=\inf_{f\in
    S(\theta)}\kappa(\Delta_f)$, where $\kappa(\Delta)$ is the
  conformal capacity of $\Delta$. Then, there is a universal constant
  $C>0$ such that $|\log(\kappa(\theta))+\Phi(\theta)|<C$.
\end{theorem}

We will also need a well-known theorem about the regularity of the
boundary of Siegel discs of quadratic polynomials. Its proof can be
found in \cite{DouadyAsterisqueSiegelHerman}.

\begin{theorem}[(Douady-Ghys)]\label{TheoremQuasiCircleDouadyGhys}
  Let $\theta$ be of bounded type, and $p(z)=e^{2\pi i
    \theta}z+z^2$. Then the boundary of the Siegel disc around 0 is a
  quasi-circle containing the critical point.
\end{theorem}
The following is a theorem by M. Herman concerning critical points on
the boundary of Siegel discs. Its proof can be found in \cite[p.
601]{HermanCriticalExponential}

\begin{theorem}[(Herman)]
  \label{TheoremHermanInjectiveBounded}
  Let $g(z)$ be an entire function such that $g(0)=0$ and
  $g'(0)=e^{2\pi i\alpha}$ with $\alpha$ Diophantine. Let $\Delta$ be
  the Siegel disc around $z=0$. If $\Delta$ has compact closure in
  $\mathbb{C}$ and $\left.g\right|_{\bar\Delta}$ is injective then $g(z)$
  has a critical point in $\partial\Delta$.
\end{theorem}

In fact, the set of Diophantine numbers could be replaced by the set
$\mathcal{H}$ of Herman numbers, where
$\mathcal{D}\subsetneq\mathcal{H}\subsetneq\mathcal{B}$, as shown in
\cite{YoccozAnalyticLinearizations}.

Finally, we state a result which is a combination of Theorems 1 and 2
in \cite{RempeSiegelAndPeriodicRays}.

\begin{definition}
  We define the class $\mathcal{B}$ as the class of entire functions
  with a bounded set of singular values.
\end{definition}

\begin{theorem}[(Rempe), \cite{RempeQuestionHBRSiegelDisks}]\label{TheoremRempeSingularitiesBoundary}
  Let $f\in\mathcal{B}$ with $S(f)\subset\curv{J}(f)$, where $S(f)$
  denotes the set of singular values of $f$. If $\Delta$ is a Siegel
  disc of $f(z)$ which is unbounded, then
  $S(f)\cap\partial\Delta\neq\emptyset$.
\end{theorem}

\subsection{Topological results}
\label{SectionPreliminarsTopologicalResults}
To prove Theorem \ref{ThmUnboundedSiegelDisks}  we need to extend a
result of Rogers in \cite{RogersSingularities} to a larger class of
functions, namely functions of finite order with no wandering domains.

The result we need follows some preliminary definitions.

\begin{definition}
  A \emph{continuum} is a compact connected non-void metric space.
\end{definition}

\begin{definition}
  A pair $(g,\Delta)$ is a local Siegel disc if $g$ is conformally
  conjugate to an irrational rotation on $\Delta$ and $g$ extends
  continuously to $\bar\Delta$.
\end{definition}

\begin{definition}
  We say a bounded local Siegel disc $(\left.f\right|_\Delta,\Delta)$
  is \emph{irreducible} if the boundary of $\Delta$ separates the
  centre of the disc from $\infty$, but no proper closed subset of the
  boundary of $\Delta$ has this property.
\end{definition}

\begin{theorem}\label{ThmRogersClassB}
  Suppose $\Delta$ is a Siegel disc of a function $f$ in the class
  $\mathcal{B}$, and $\partial \Delta$ is a decomposable continuum.
  Then $\partial \Delta$ separates $\mathbb{C}$ into exactly two
  complementary domains.
\end{theorem}


For the proof of this theorem we will need the following ingredients
which will be only used in this proof. The topological results can be
found in any standard reference on algebraic topology.

\begin{theorem}\label{ThmRogersMainTheorem}
  If $(\Delta,f_\theta)$ is a bounded irreducible local Siegel disc, then
  the following are equivalent:
\begin{itemize}
\item $\partial \Delta$ is a decomposable continuum,
\item each pair of impressions is disjoint, and
\item the inverse of the map $\varphi:\mathbb{D}\to \Delta$ extends
  continuously to a map $\Psi:\partial$ $\Delta\to S^1$ such that for
  each $\eta\in S^1$, the fibre $\Psi^{-1}(\eta)$ is the impression
  $I(\eta)$.
\end{itemize}
\end{theorem}
\begin{proof}
  See \cite{RogersSingularities}.

\end{proof}

\begin{theorem}[(Vietoris-Begle)]
  Let $X$ and $Y$ be compact metric spaces and $f:X\to Y$ continuous and
  surjective and suppose that the fibres are acyclic, i.e.
  \[\tilde H^r(f^{-1}(y))=0, 0\leq r\leq n-1, \quad\forall y\in Y,\]
  where $\tilde H^r$ denotes the $r$-th reduced co-homology group.
  Then, the induced homomorphism
  \[f^*:\tilde H^r(Y)\to \tilde H^r(X)\] is an isomorphism
  for $r\leq n-1$ and is a surjection for $r=n$.
\end{theorem}

\begin{theorem}[(Alexander's duality)]
  Let $X$ be a compact sub-space of the Euclidean space $E$ of
  dimension $n$, and $Y$ its complement in $E$. Then,
  \[\tilde H_q(X)\cong \tilde H^{n-q-1}(Y)\] where $\tilde
  H_*,\,\tilde H^*$ stands for \v Cech reduced homology and reduced
  co-homology respectively.
\end{theorem}

\begin{remark}\label{RemarkAlexander}
  The case $E=S^2$, $X=S^1$ (or $H^1(X)=\mathbb{Z}$) is Jordan's Curve
  Theorem.
\end{remark}

\begin{definition}
  If $X$ is a compact subset of $\mathbb{C}$, then the three following
  conditions are equivalent:
\begin{itemize}
\item $X$ is cellular,
\item $X$ is a continuum that does not separate $\mathbb{C}$,
\item $H^1(X)=0=\tilde H^0(X)$,
\end{itemize}
where $\tilde H^r(X)$ stands for reduced \v Cech co-homology and
$H^r(X)$ for \v Cech co-homology.
\end{definition}

\begin{definition}
  We say a map $f:X\to Y$ is cellular if each fibre $f^{-1}(y)$ is a
  cellular set.
\end{definition}

\begin{remark}
  Recall that $\tilde H^1(X)\cong H^1(X)$.
\end{remark}

\begin{remark}\label{RemarkVietorisCellular}
  By definition and in view of the Vietoris-Begle Theorem, cellular
  maps induce isomorphisms between first reduced co-homology groups.
\end{remark}

\paragraph*{\it \normalsize Proof of Theorem \ref{ThmRogersClassB}.}
We first show that any Siegel disc $\Delta$ for $f\in\mathcal{B}$ is a
bounded irreducible local Siegel disc. Recall that we define the
escaping set of a function $f:\mathbb{C}\to\mathbb{C}$ as:
\[I(f)=\{z\vert\, f^n(z)\to\infty\textrm{ as }n\to\infty\}.\] Clearly
$(\left.f\right|_\Delta,\Delta)$ is a local Siegel disc. It is also
bounded by assumption. The only thing left to prove is it is
irreducible. If $X$ is a proper closed subset of $\partial\Delta$ and
if $x$ is a point of $\partial\Delta\backslash X$, then there is a
small disc $B$ containing $x$ and missing $X$. Since
$x\in\partial\Delta$, the disc $B$ contains a point of $\Delta$. As
$x\in\partial \Delta\subset \curv{J}(f)$, the disc $B$ contains a
point $y\in I(f)$. Now, Theorem 3.1.1 in \cite{RottenfusserThesis}
states that for $f\in\mathcal{B}$ the set $I(f)\cup\{\infty\}$ is
arc-connected, and thus $y$ can be arc-connected to $\infty$ through
points in $I(f)$. It follows that the centre of the Siegel disc and
infinity are in the same complementary domain of $\mathbb{C}\backslash
X$.

Clearly $\Psi(\eta)$ for $\eta\in S^1$ is a continuum, which is called
the impression of $\eta$ and denoted $\mathrm{Imp}(\eta).$ Furthermore,
$\mathrm{Imp}(\eta)$ does not separate $\mathbb{C}$. Indeed, if $U$ is a
bounded complementary domain of $\mathrm{Imp}(\eta)$, then either $f^n(U)\cap
U=\emptyset$ for all $n$ or there are intersection points. Clearly
$f^n(U)\cap U=\emptyset$, as if $f^n(U)\cap U\neq\emptyset$ for some
$n$, then $f^n(\partial U)\cap \partial U\neq\emptyset$, but this
implies $\mathrm{Imp}(\eta)=F^n(\mathrm{Imp}(\eta))=\mathrm{Imp}(\eta+n\theta)$ and as $\partial
\Delta$ is a decomposable continuum, each pair of impressions is
disjoint by Theorem \ref{ThmRogersMainTheorem} and this intersection
must be empty. Hence, $f^n(U)\cap U=\emptyset$ for all
$n\in\mathbb{N}$ which implies $U$ is a wandering domain, and for
functions in $\mathcal{B}$ it is known there are no wandering domains
(see \cite{EremenkoLyubichDynamicalEntire92}).

Therefore $\mathrm{Imp}(\eta)$ is a cellular set and thus $\Psi$ is a cellular
map. The Vietoris-Begle theorem implies that the induced homomorphism
$\Psi^*:\tilde H^1(S^1)\to \tilde H^1(\partial \Delta)$ is an
isomorphism (see Remark \ref{RemarkVietorisCellular}). Then $\tilde
H^1(\partial \Delta)=\mathbb{Z}$ and by Alexander's duality $\partial
\Delta$ separates $\mathbb{C}$ into exactly two complementary domains
(see Remark \ref{RemarkAlexander}).
  
  \hfill
  \qed



\section{The (entire transcendental) family $f_a$}
\label{SectionCharacterizationImagesDynamicalPlanes}
In this section we describe the dynamical plane of the family of
entire transcendental maps
\[
f_a(z) = \lambda a(e^{z/a}(z+1-a)-1+a),
\]
for different values of $a\in \mathbb{C}^*$, and for $\lambda=e^{2\pi
  i\theta}$, with $\theta$ being a fixed irrational Brjuno number
(unless otherwise specified). For these values of $\lambda$, in view
of Theorem \ref{TheoremLinearizationOfBrjunoRussmann} there exists an
invariant Siegel disc around $z=0$, for any value of $a \in
\mathbb{C}^*$.

We start by showing that this family contains all possible entire
transcendental maps with the properties we require.

\begin{theorem}
\label{TheoremCharacterizationOfTheFamily}
Let $g(z)$ be an entire transcendental function having the following
properties
\begin{enumerate}
\item finite order,
\item one asymptotic value $v$, with exactly one finite pre-image $p$
  of $v$,
\item a fixed point (normalised to be at 0) of multiplier $\lambda \in
  \mathbb{C}$,\label{ItemCharacterizationCondition0Fixed}
\item a simple critical point (normalised to be at $z=-1$) and no
  other critical points.
\end{enumerate} 
Then $g(z)=f_a(z)$ for some $a\in\mathbb{C}$ with $v=\lambda a (a-1)$ and
$p=a-1$. Moreover no two members of this family are conformally
conjugate.
\end{theorem}

\begin{proof}
  As $g(z)-v=0$ has one solution at $z=p$, we can write:
  \[g(z)=(z-p)^m e^{h(z)}+v,\] where, by Hadamard's factorisation
  theorem (Theorem \ref{TheoremHadamardsFactorization}), $h(z)$ must
  be a polynomial, as $g(z)$ has finite order. The derivative of this
  function is

\[g'(z)=e^{h(z)}(z-p)^{m-1}(m+(z-p)h'(z)),\] whose zeroes are the
solutions of $z-p=0$ (if $m>1$) and the solutions of
$m+(z-p)h'(z)=0$. But as the critical point must be simple and unique,
$m=1$ and $\deg h'(z)=0$.
Therefore
\[g(z)=(z-p)e^{\alpha z+\beta}+v,\] and from the expression for the
critical points,
\[
\alpha=\frac{1}{p+1}.
\]
Moreover from the fact that $g(0)=0$ we can deduce that $v=pe^\beta$,
and from condition \ref{ItemCharacterizationCondition0Fixed},
i.e. $g'(0)=\lambda$, we obtain $e^\beta=\lambda(1+p)$.  All together
yields
\[g(z)=\lambda (z-p)(1+p)e^{z/(1+p)}+\lambda p(1+p).\]
Writing $a=p+1$ we arrive to
\[g(z)=\lambda a(z-a+1)e^{z/a}+\lambda a(a-1)=f_a(z),\]
as we wanted.

Finally, if $f_a(z)$ and $f_{a'}(z)$ are conformally conjugate, the
conjugacy must fix 0,-1 and $\infty$ and therefore is the identity
map.

\end{proof}


\subsection{Dynamical planes}

For any parameter value $a\in\mathbb{C}^*$, the Fatou set always
contains the Siegel disc $\Delta_a$ and all its pre-images. Moreover,
one of the singular orbits must be accumulating on the boundary of
$\Delta_a$. The other singular orbit may then either eventually fall
in $\Delta_a$, or accumulate in $\partial\Delta_a$, or have some
independent behaviour. In the first case we
say that the singular value is captured by the Siegel disc. More
precisely we define the \emph{capture parameters} as
\begin{align}
  C=\{a\in\mathbb{C}^*\vert &f_a^n(-1)\in\Delta_a\textrm{ for some
  }n\geq 1 \textrm{ or } \nonumber\\
  &f_a^n(\asv{a})\in\Delta_a\textrm{ for some }n\geq 0\}\nonumber
\end{align}
Naturally $C$ splits into two sets $C=C^c\cup C^v$ depending on
whether the captured orbit is the critical orbit ($C^c$) or the orbit
of the asymptotic value ($C^v$). We will follow this convention,
superscript $c$ for critical and superscript $v$ for asymptotic,
throughout this paper.

In the second case, that is, when the free singular value has an
independent behaviour, it may happen that it is attracted to an
attracting periodic orbit. We define the \emph{semi-hyperbolic
  parameters} $H$ as
\[H=\{a\in\mathbb{C}^*\vert f_a\textrm{ has an attracting periodic
  orbit}\}.\] Again this set splits into two sets, $H=H^c\cup H^v$
depending on whether the basin contains the critical point or the
asymptotic value.

Notice that these four sets $C^c,\,C^v,\,H^c,\,H^v$ are pairwise
disjoint, since a singular value must always belong to the Julia set,
as its orbit has to accumulate on the boundary of the Siegel disc.

In the following sections we will describe in detail these regions of
parameter space, but let us first show some numerical experiments. For
all figures we have chosen $\theta=\frac{1+\sqrt{5}}{2}$, the golden
mean number.

Figure \ref{FigureParametersEscapeToInfinityZoomOut} (in the
Introduction) shows the parameter plane, where the left side is made
with a simple escaping algorithm. The component containing $a=1$ is
the main capture component for which $\asv{a}$ itself belongs to the
Siegel disc. On the right side we see the same parameters, drawn with
a different algorithm. Also in Figure
\ref{FigureParametersEscapeToInfinityZoomOut}, we can partially see
the sets $H^v_1$ and $H^v_2$ (and infinitely many others), where the
sub-indices denote the period of the attracting orbit.

\begin{figure}[!hbt]
  \begin{center}
    \fbox{\includegraphics[width=0.45\textwidth]{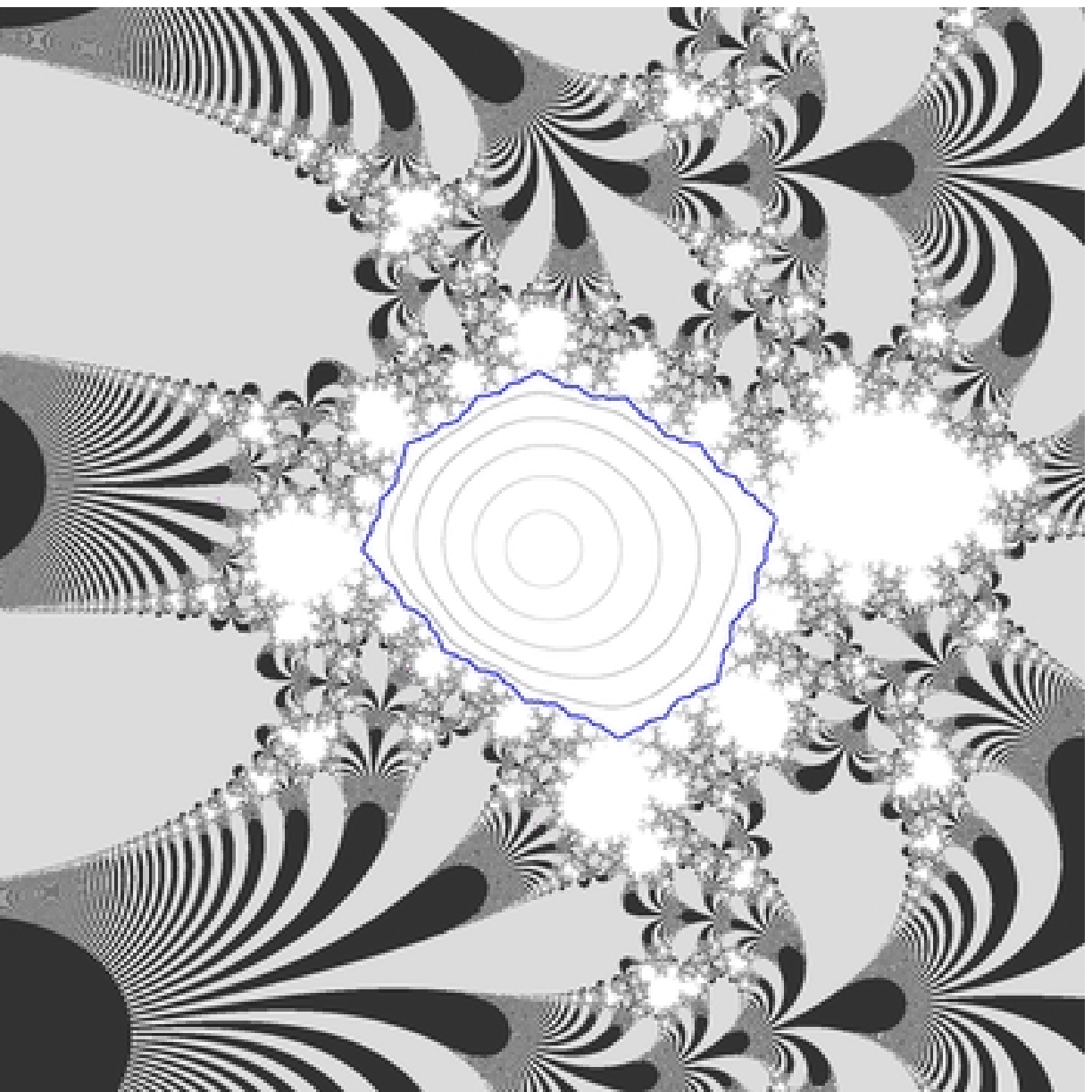}
    }
    \fbox{\includegraphics[width=0.45\textwidth]{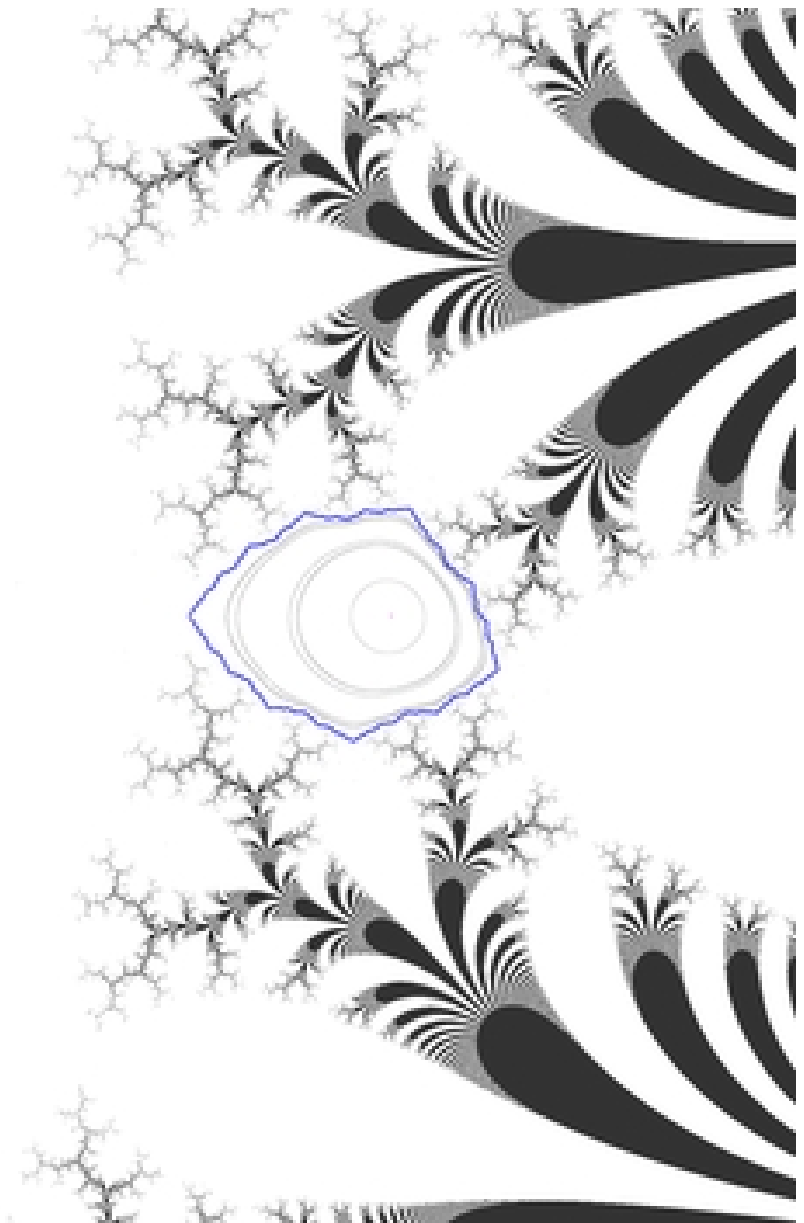} }
    \caption{{\bf Left:} Julia set for a parameter in a
      semi-hyperbolic component (for the asymptotic value). Details:
      $a=(-0.62099,0.0100973),$ upper left: $(-4,3),$ lower right:
      $(2,-3)$. In light grey we see the attracting basin of the
      attracting cycle, and in white the Siegel disc and its
      pre-images. {\bf Right:} Julia set of the semi-standard map,
      corresponding to $f_1(z)=\lambda ze^z$. Upper left: $(-3,3),$
      lower right: $(3,-3)$. The boundary of the Siegel disc around 0
      is shown, together with some of the invariant curves. The Fatou
      set consists exclusively of the Siegel disc and its
      pre-images.}
    \label{FigureJuliaSetHyperbolicNotPolyLikeAndSemiStandard}
  \end{center}
\end{figure}

In Figure \ref{FigureJuliaSetHyperbolicNotPolyLikeAndSemiStandard}
(left) we can see the dynamical plane for $a$ chosen in one of the
semi-hyperbolic components of Figure
\ref{FigureParametersEscapeToInfinityZoomOut}, where the Siegel disc
and the attracting orbit and corresponding basin are shown in
different colours.

Figure \ref{FigureJuliaSetHyperbolicNotPolyLikeAndSemiStandard}
(right) shows the dynamical plane of $f_1(z)=\lambda ze^{z}$, the
semi-standard map. In this case the asymptotic value $\asv{1}=0$ is
actually the centre of the Siegel disc. It is still an open question
whether, for some exotic rotation number, this Siegel disc can be
unbounded. For bounded type rotation numbers, as the one in the
figure, the boundary is a quasi-circle and contains the critical point
\cite{GeyerSiegelHerman}.

\begin{figure}[!hbt]
  \begin{center}
    \psfrag{0}{$0$}
    \fbox{\includegraphics[width=0.45\textwidth]{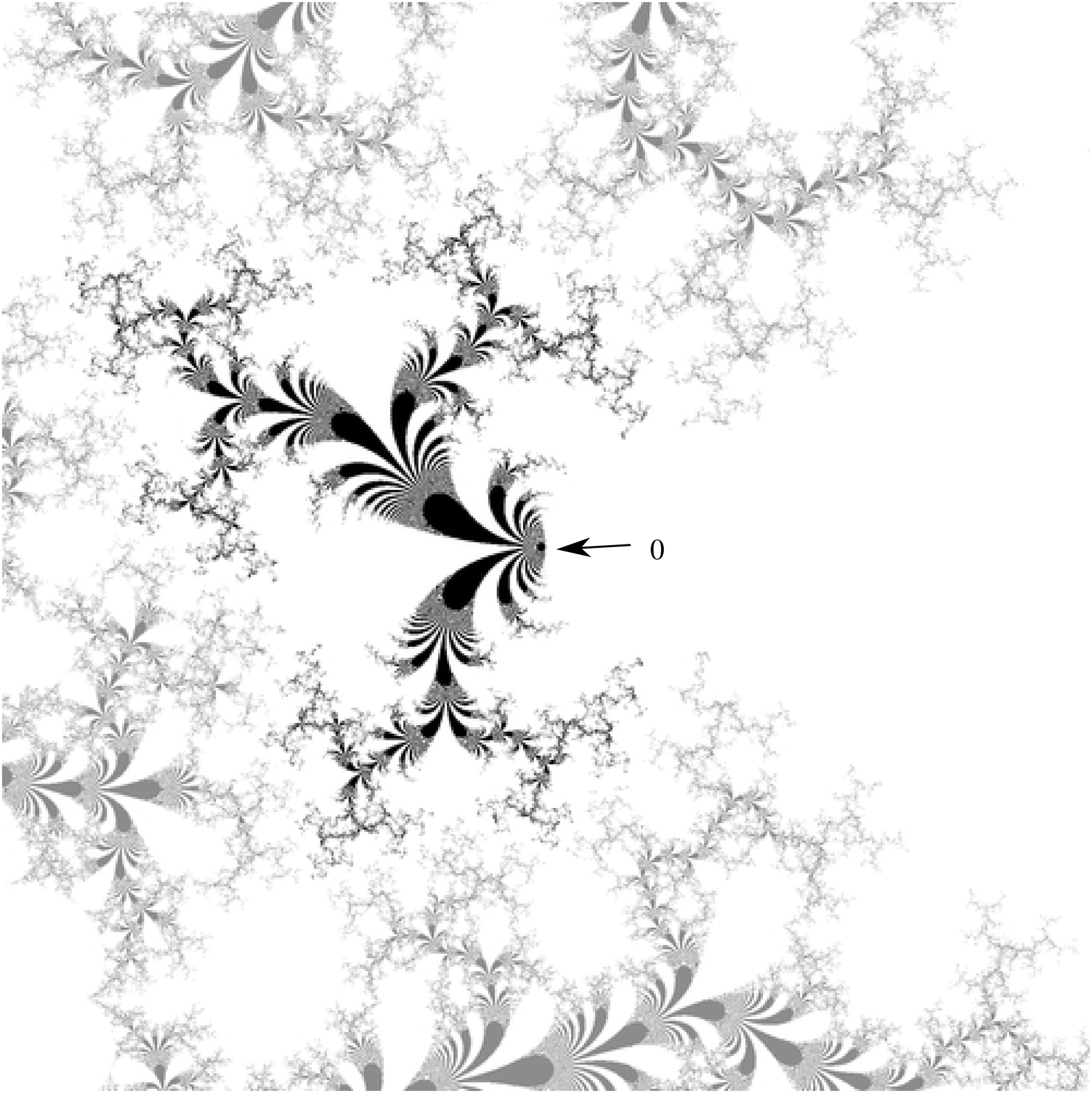} }
    \fbox{\includegraphics[width=0.45\textwidth]{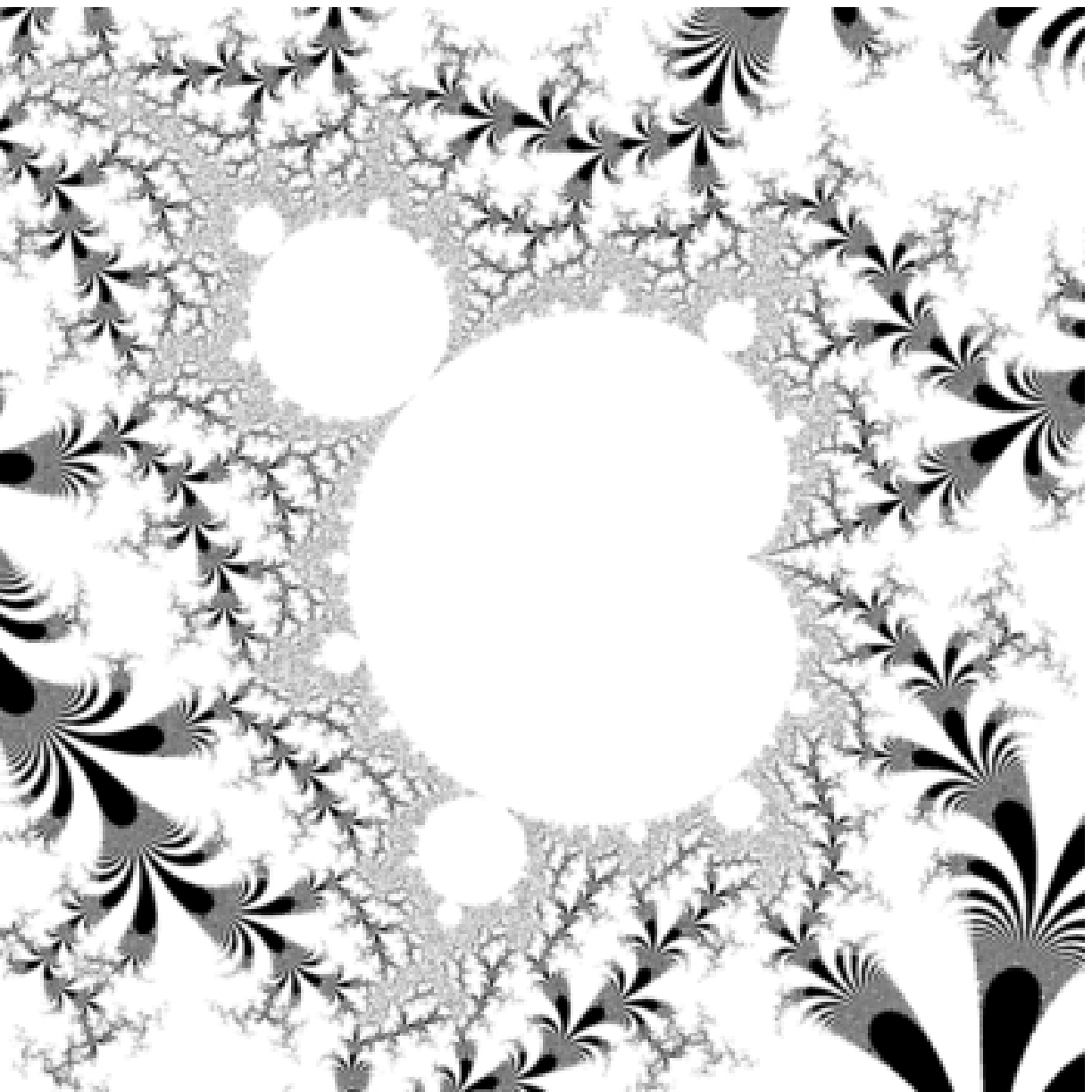}
    }
    
    \label{FigureParameterPlaneZoomBabyMandelbrot}
    \caption{ {\bf Left:} ``Crab''-like structure corresponding to
      escaping critical orbits (dark grey). Upper left: $(-0.6,0.6),$
      lower right: $(0.6,-0.6)$. In light grey we see parameters for
      which the orbit of $\asv{a}$ escapes. {\bf Right:} Baby
      Mandelbrot set from a close-up in the ``crab like'' structure.
      Upper left: $(-0.336933,0.1128),$ lower right:
      $(-0.322933,0.08828)$.}
    \label{FigureParametersEscapeToInfinityZoomCrab}
  \end{center}
\end{figure}

Figure \ref{FigureParametersEscapeToInfinityZoomCrab}, left side,
shows a close-up view of the parameter region around $a=0$, and in the
right side, we can see a closer view of a random spot, in
particular a region in $H^c$, that is, parameters for which the
critical orbit is attracted to a cycle.

\begin{figure}[hbt!]
  \begin{center}
    \fbox{\includegraphics[width=0.45\textwidth]{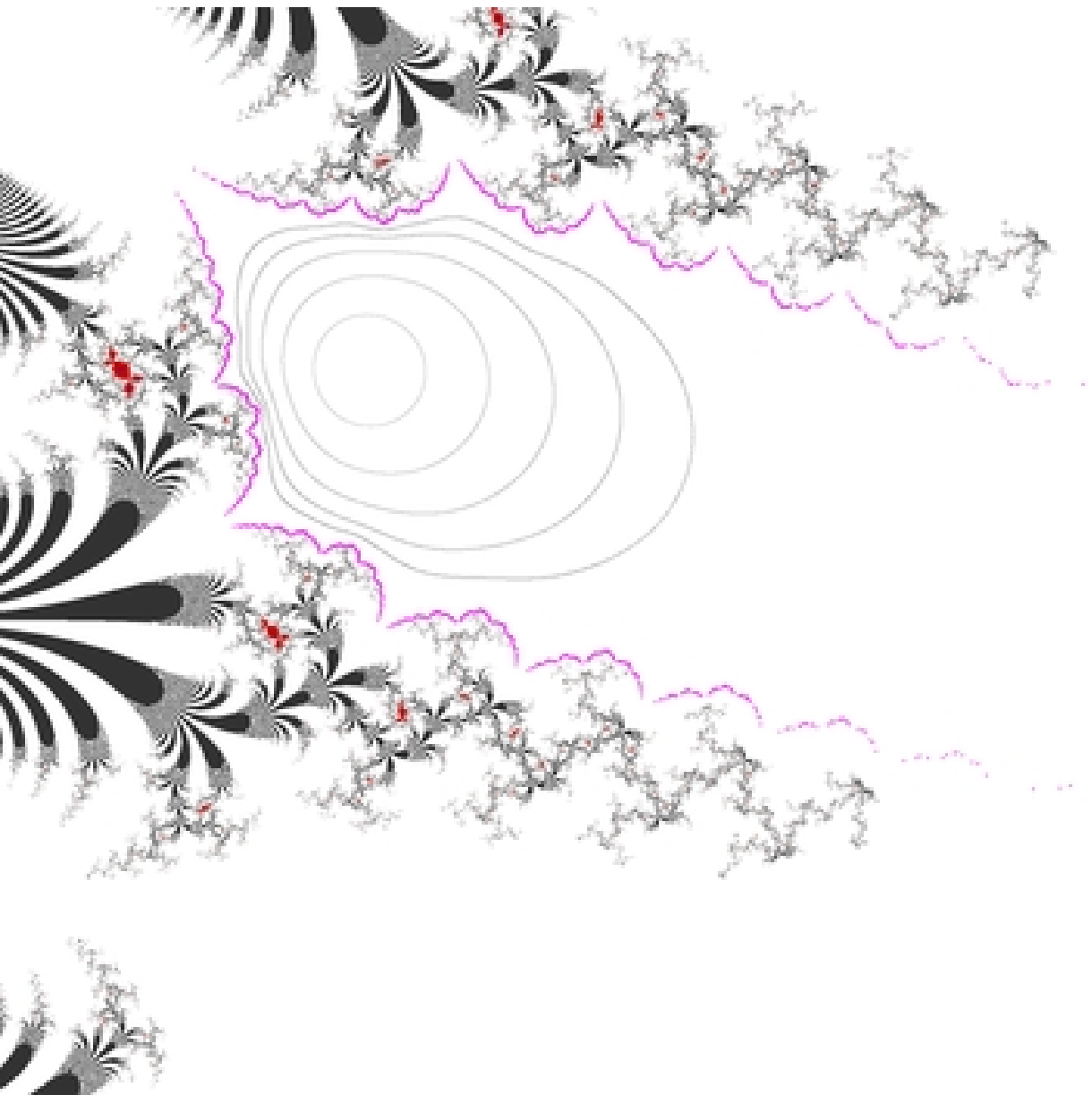}}
    \fbox{\includegraphics[width=0.45\textwidth]{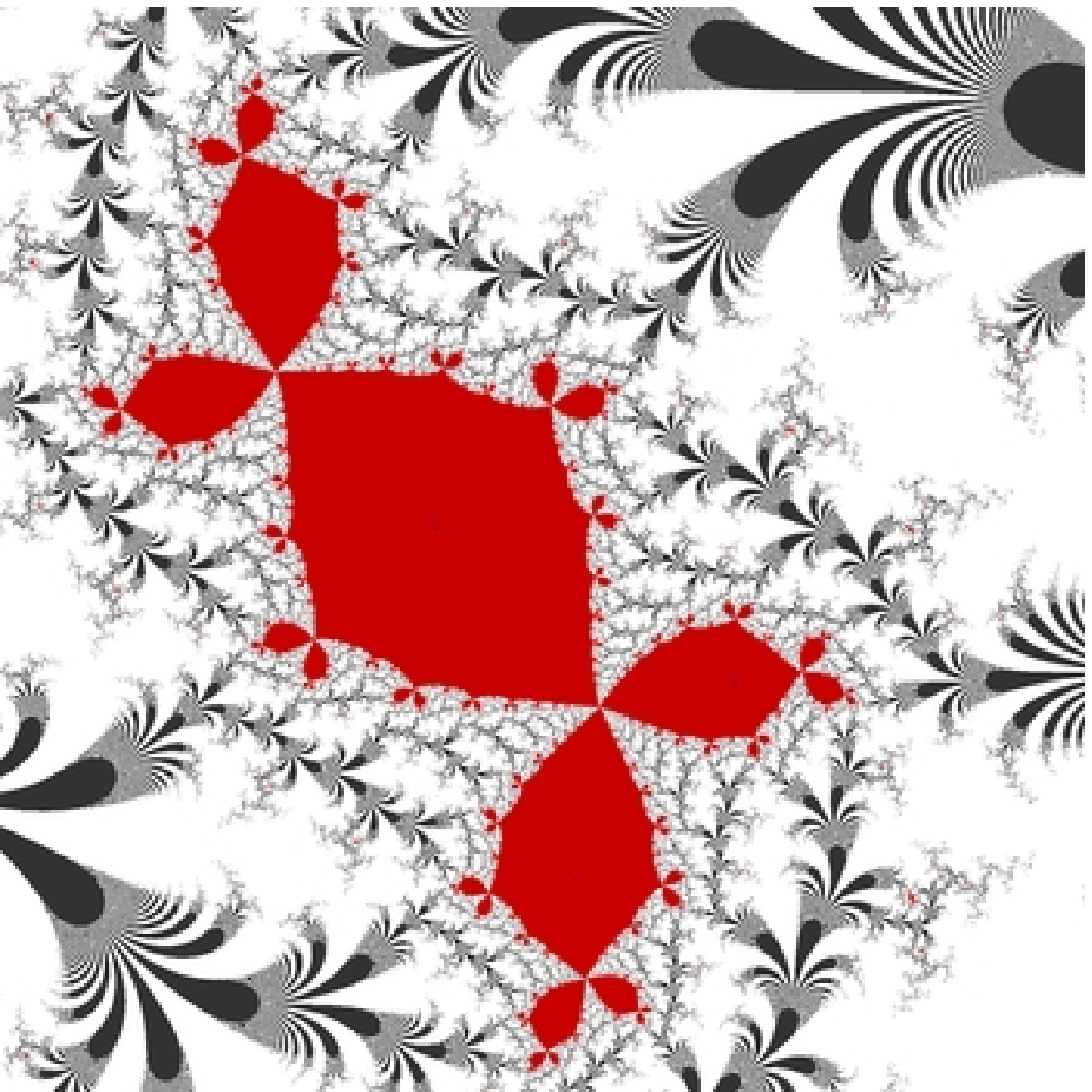}}
    \caption{ {\bf Left:} Julia set for a parameter in a
      semi-hyperbolic component for the critical value. By Theorem
      \ref{ThmUnboundedSiegelDisks} this Siegel disc is unbounded.
      Details: $a=(-0.330897,0.101867),$ upper left: $(-1.5, 1.5).$,
      lower right: $[3,-3]$. {\bf Right:} Close-up of a basin of
      attraction of the attracting periodic orbit. Upper left:
      $(-1.1,0.12),$ lower right: $(-0.85,-0.13)$.}
    \label{FigureJuliaSetUnboundedSiegelDiskCenterMandelbrot}
  \end{center}
\end{figure}

One of these dynamical planes is shown in Figure
\ref{FigureJuliaSetUnboundedSiegelDiskCenterMandelbrot}. Observe that
the orbit of the asymptotic value is now accumulating on
$\partial\Delta_a$ and we may have unbounded Siegel discs.

\begin{figure}[!hbt]
  \begin{center}
    \fbox{\includegraphics[width=0.45\textwidth]{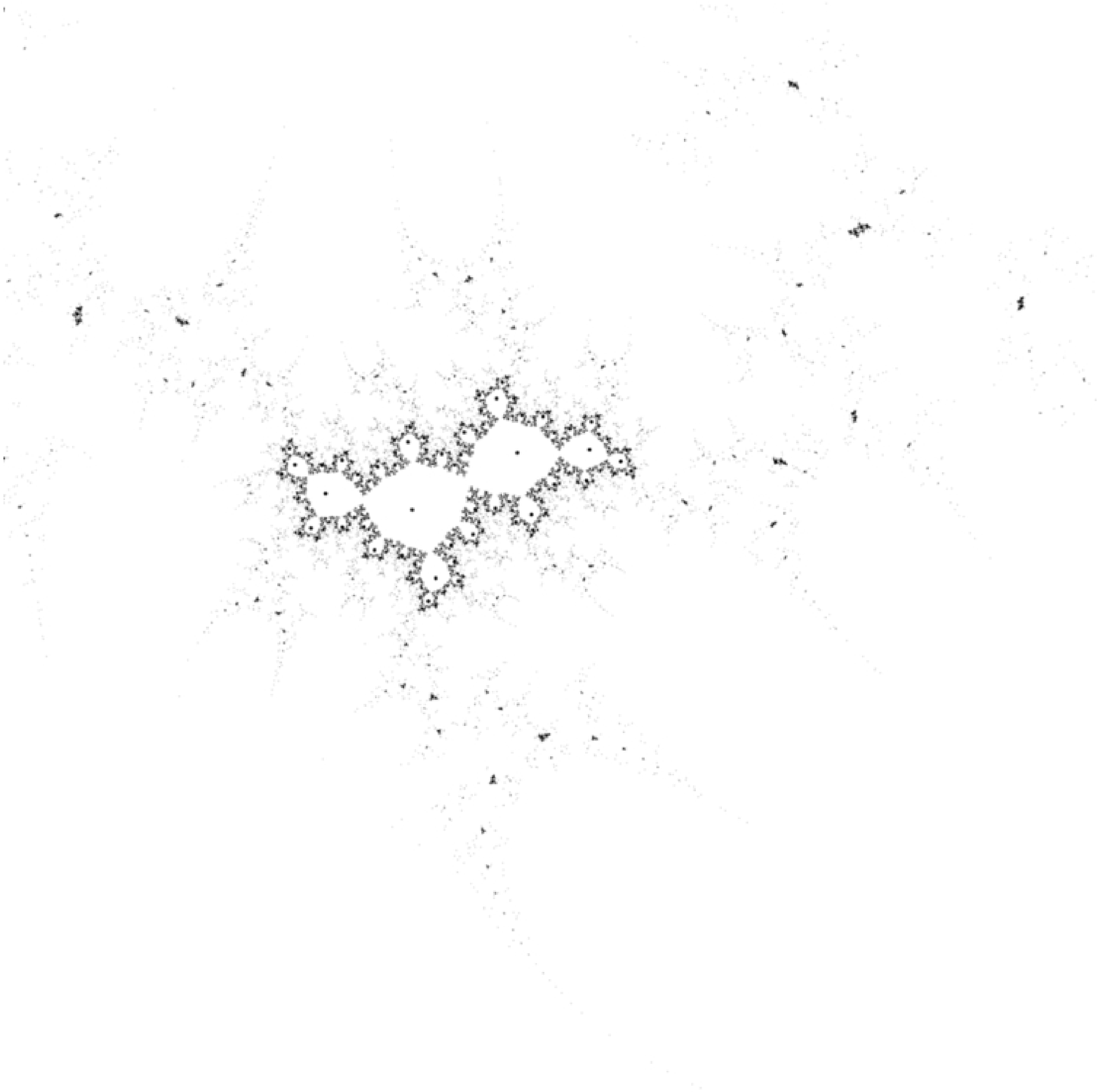}}
    \caption{A close up of Figure
      \ref{FigureParametersEscapeToInfinityZoomOut}, Right. A
      quadratic Siegel disc in parameter space, corresponding to a
      capture zone for the asymptotic value. Upper left:
      $(7.477,4.098),$ Lower right: $(7.777,3.798)$.}
  \label{FigureParametersSiegelDisk}
  \end{center}
\end{figure}

Finally Figure \ref{FigureParametersSiegelDisk} shows some components
of $C^v$, where the orbit of the asymptotic value is captured by the
Siegel disc.


We start by considering large values of $a\in\mathbb{C}^*$. By
expanding $f_a(z)$ into a power series it is easy to check that as
$a\to\infty$ the function approaches the quadratic polynomial $\lambda
z(1+z/2)$. It is therefore not surprising that we have the following
theorem, which we shall prove at the end of this section.

\begin{theorem}\label{TheoremTheFamilyIsPolynomialLike}
  There exists $M>0$ such that the entire transcendental family
  $f_a(z)$ is polynomial-like of degree two for $|a|>M$. Moreover, the
  Siegel disc $\Delta_a$ (and in fact, the full small filled Julia
  set) is contained in a disc of radius $R$ where $R$ is a constant
  independent of $a$.
\end{theorem}

\begin{figure}[!hbt]
  \begin{center}
    \fbox{\includegraphics[width=0.45\textwidth]{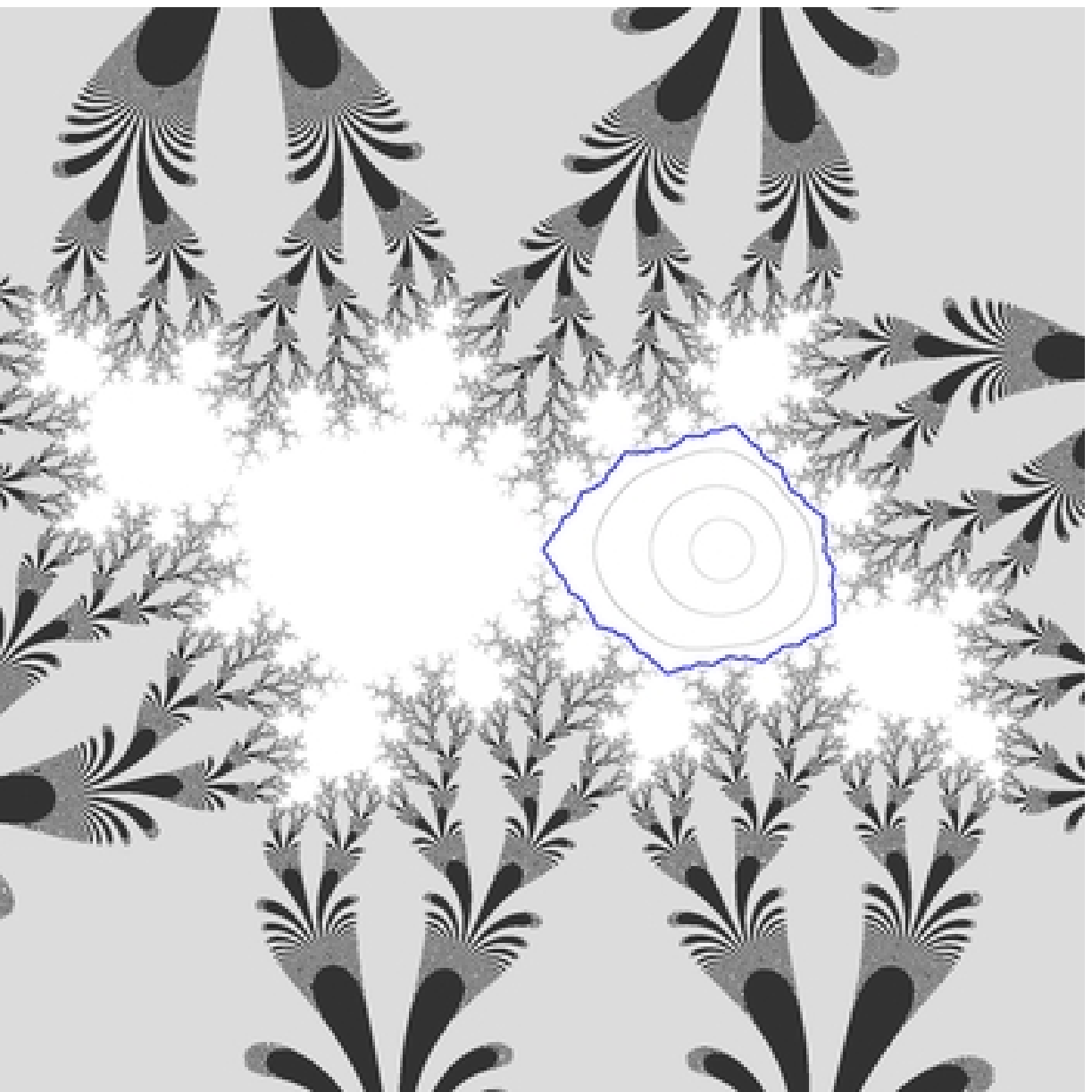}}
    \fbox{\includegraphics[width=0.45\textwidth]{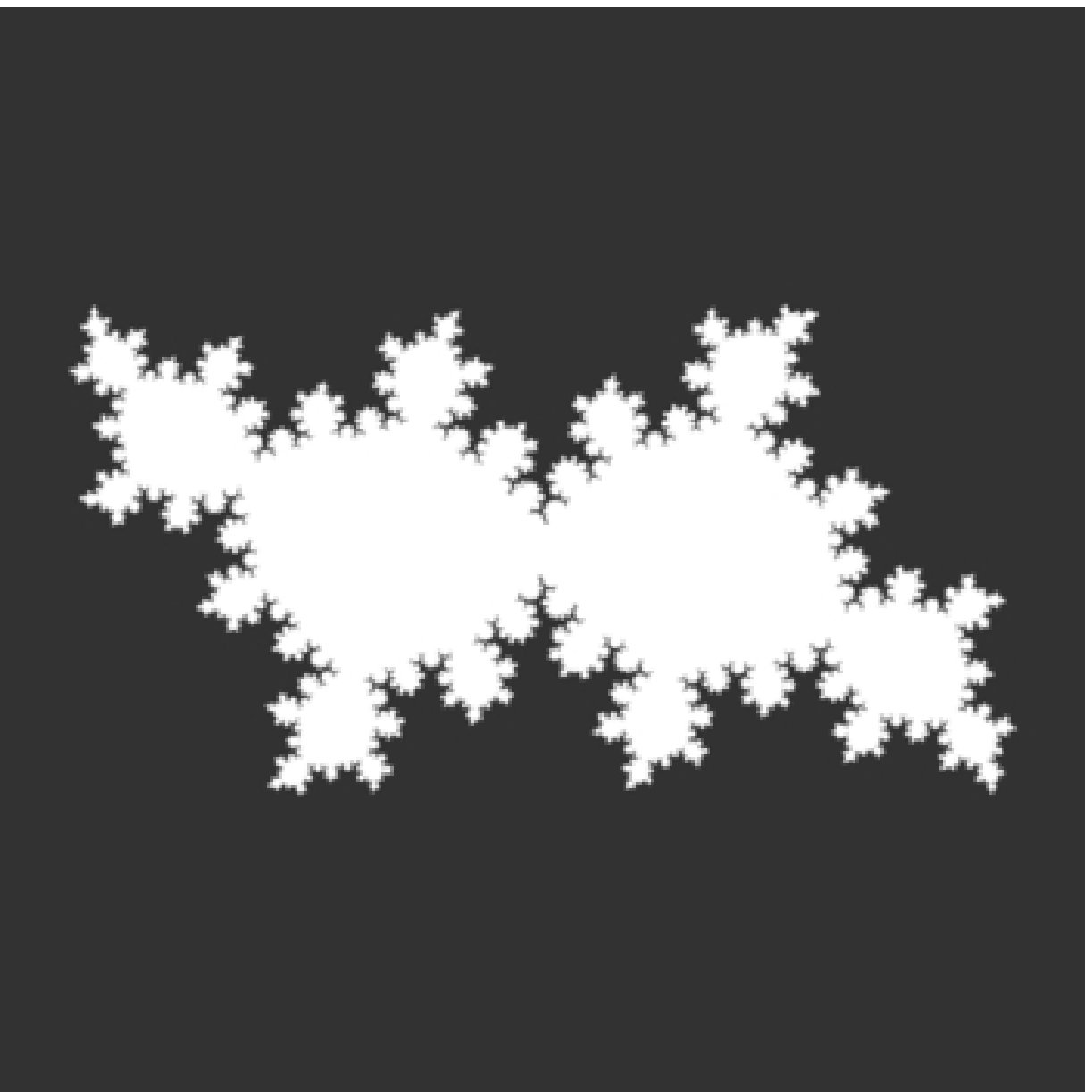}}
    \caption{ {\bf Left:} Julia set corresponding to a polynomial-like
      mapping. Details: $a=(15,-15),$ upper left: $(-4,3),$ lower
      right: $(2,-3)$. {\bf Right:} Julia set corresponding to the related
      polynomial. Upper left: $(-4,3),$ lower right: $(-2,3)$}
    \label{FigureJuliaSetForPolynomialLike}
  \end{center}
\end{figure}

Figure \ref{FigureJuliaSetForPolynomialLike} shows the dynamical plane
for $a=15+15i$, $\lambda=e^{2\pi(\frac{1+\sqrt{5}}{2})i}$ where we clearly see
the Julia set of the quadratic polynomial $\lambda z(1+z/2)$, shown on
the right side.

An immediate consequence of Theorem
\ref{TheoremTheFamilyIsPolynomialLike} above follows from Theorem
\ref{TheoremQuasiCircleDouadyGhys}. This is Part \ref{PartATmaA}) of
Theorem A in the Introduction.

\begin{corollary}\label{CorollaryConstantTypeQuasiCircle}
  For $|a|>M$, and $\theta$ of constant type the boundary of
  $\Delta_a$ is a quasi-circle that contains the critical point.
\end{corollary}

In fact we will prove in Section
\ref{SectionPropertiesCaptureComponents} (Proposition
\ref{PropositionHolomorphicMotionOfC0}) that the same occurs in many
other situations like, for example, when the asymptotic value lies
itself inside the Siegel disc or when it is attracted to an attracting
periodic orbit. See Figures
\ref{FigureJuliaSetHyperbolicNotPolyLikeAndSemiStandard} (Left) and
\ref{FigureJuliaSetForPolynomialLike}.

\begin{figure}[!hbt]
  \begin{center}
    \fbox{ \includegraphics[width=0.45\textwidth]{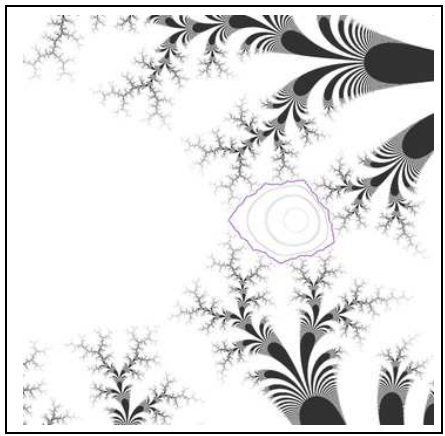}}
    \caption{ Julia set for the parameter
      $a\approx1.544913893+0.32322773i$. The parameter is chosen so
      that the critical point and the asymptotic value are at the same
      point, hence both singular orbits accumulate on the boundary.
      Upper left: $(-1.5,1.5),$ lower right: $(3,-3)$.}
    \label{FigureSiegelDiskBoundedWithAsymptoticValueInTheBoundary}
  \end{center}
\end{figure}

In fact we believe that this family provides examples of Siegel discs
with an asymptotic value on the boundary, but such that the boundary
is a quasi-circle containing also the critical point. A parameter
value with this property could be given by
$a_0\approx1.544913893+0.32322773i\in\partial C_0^v$,
$\lambda=e^{2\pi(\frac{1+\sqrt{5}}{2})i}$ (see Figure
\ref{FigureSiegelDiskBoundedWithAsymptoticValueInTheBoundary}) where
the asymptotic value and the critical point coincide.

\begin{figure}[hbt!]
  \begin{center}
    \fbox{\includegraphics[width=0.45\textwidth]{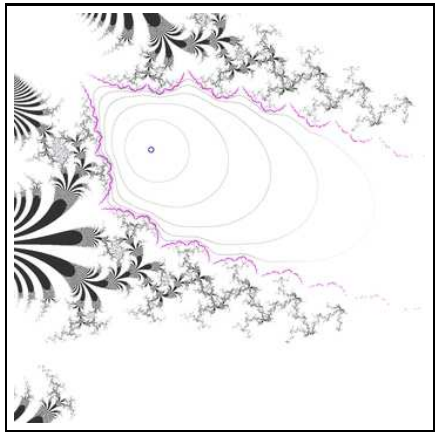}}
    \caption{ Point in a capture component for the critical value, so
      that the Siegel disc is either unbounded or an indecomposable
      continuum. Details: $a=(-0.33258,0.10324),$ upper left: $(-1.5,
      1.5)$, lower right: $(-3,-3)$.}
    \label{FigureJuliaSetUnboundedSiegelDiskCenterCaptureCritical}
  \end{center}
\end{figure}

The opposite case, that is, the Siegel disc being unbounded and its
boundary non-locally connected also takes place for certain values of
the parameter $a$, as we show in the following theorem, which covers
parts \ref{PartBTmaA}) and \ref{PartCTmaA}) of Theorem A (see Figure
\ref{FigureJuliaSetUnboundedSiegelDiskCenterCaptureCritical}).

\begin{theorem}
  \label{ThmUnboundedSiegelDisks}
  
  Let $\theta$ be Diophantine\footnote{Diophantine numbers can
    actually be replaced by the larger class of irrational numbers
    $\mathcal{H}$ (see \cite{YoccozAnalyticLinearizations},
    \cite{PerezMarcoFixedpoints})}, then:
  \begin{enumerate}[a)]
  \item \label{ThmUnbddSiegelA} If $f_a^n(-1)\to\infty$ then
    $\Delta_a$ is unbounded and $\asv{a}\in\partial\Delta_a$,
  \item \label{ThmUnbddSiegelB} if $a\in H^c\cup C^c$ either $\Delta_a$
    is unbounded or $\partial\Delta_a$ is an indecomposable continuum.
  \end{enumerate}
\end{theorem}

\begin{proof}
  The proof of the first part is a slight modification of Herman's
  proof for the exponential map (see
  \cite{HermanCriticalExponential}). The difference is given by the
  fact that the asymptotic value of $f_a(z)$ is not an omitted value,
  and by the existence of a second singular value. For both parts we
  need the following definitions. Suppose that $\Delta:=\Delta_a$ is
  bounded and let $\Delta_i$ denote the bounded components of
  $\mathbb{C}\backslash\partial\Delta$. Let $\Delta_\infty$ be the
  unbounded component. Since $\Delta$ and $\Delta_i$ are simply
  connected, then $\hat\Delta:=\mathbb{C}\backslash\Delta_\infty$ is
  compact and simply connected. By the Maximum Modulus Principle and
  Montel's theorem,
  $\{\left.f_a^n\right|_{\Delta_i}\}_{n\in\mathbb{N}}$ form a normal
  family and hence $\Delta_i$ is a Fatou component. We also have that
  $\partial\Delta=\partial\Delta_\infty$, although this does not imply
  a priori that $\Delta_i=\emptyset$ (see Wada lakes and similar
  examples \cite{RogersWadaLakes}).
  \paragraph*{\normalsize \it Proof of Part \ref{ThmUnbddSiegelA}).}
  Now suppose the critical orbit is unbounded. Then
  $c\in\mathcal{J}(f_a)$, but $\hat \Delta\cap\mathcal{J}(f_a)$ is
  bounded and invariant. Hence $c\notin\hat\Delta$.

  We claim that there exists $U$ a simply connected neighbourhood of
  $\hat\Delta$ such that $U$ contains no singular values. Indeed,
  suppose that the asymptotic value $\asv{a}$ belongs to $\hat
  \Delta$. Since $\asv{a}\in\mathcal{J}(f)$, then
  $\asv{a}\in\partial\Delta$. But $\Delta$ is bounded, and
  $\left.f\right|_{\partial\Delta}$ is surjective, hence the only
  finite pre-image of $\asv{a}$, namely $a-1$, also belongs to
  $\partial\Delta$. This means that $\asv{a}$ is not acting as an
  asymptotic value but as a regular point, since $f(z)$ is a local
  homeomorphism from $a-1$ to $\asv{a}$.

  Hence there are no singular values in $U$. It follows that
  \[\left.f\right|_{f^{-1}(U)}:f^{-1}(U)\to U\] is a covering and
  $f^{-1}:\Delta\to\Delta$ extends to a continuous map $h(z)$ from
  $\bar\Delta$ to $\bar\Delta$. Since $hf=fh=id$, it follows that
  $\left.f\right|_{\partial\Delta}$ is injective. As this mapping is
  always surjective, it is a homeomorphism. We now apply Herman's main
  theorem in \cite{HermanCriticalExponential} (see Theorem
  \ref{TheoremHermanInjectiveBounded}) to conclude that
  $\partial\Delta$ must have a critical point, which contradicts our
  assumptions. It follows that $\Delta$ is unbounded. Finally Theorem
  \ref{TheoremRempeSingularitiesBoundary} implies that
  $\asv{a}\in\partial\Delta_a$.
  \paragraph*{\it\normalsize Proof of part \ref{ThmUnbddSiegelB}).} 
  The work was done already when proving Theorem
  \ref{ThmRogersClassB}. Since $f_a$ has 2 singular values, it belongs
  to the Eremenko-Lyubich class $\mathcal{B}$. Hence, if we assume
  that $\Delta_a$ is bounded, it follows from Theorem
  \ref{ThmRogersClassB} that $\partial\Delta_a$ is either an
  indecomposable continuum or $\partial\Delta_a$ separates
  $\mathbb{\hat C}$ in exactly two complementary domains. This would
  imply that $\hat\Delta=\bar\Delta$ and by hypothesis
  $-1\notin\bar\Delta$. The same arguments as in Part
  \ref{ThmUnbddSiegelA} concludes the proof.
\end{proof}

\begin{remark}
  In part \ref{ThmUnbddSiegelA}) it is not strictly necessary that the
  critical orbit tends to infinity. In fact we only use that the
  critical point is in $\curv{J}(f_a)$ and some element of its orbit
  belongs to $\Delta_\infty$.
\end{remark}

\subsection{Large values of $|a|$: Proof of theorem
  \ref{TheoremTheFamilyIsPolynomialLike}}

Let $D:=\{w\in\mathbb{C}\vert |w|<R\},$ $\gamma=\partial D,$
$g(z)=\lambda z(z/2+1).$ If we are able to find some $R$ and $S$
such that
\begin{align}
  \left|g(z)-w\right|_{\substack{z\in\gamma\\w\in D}}&\geq S,\nonumber\\
  \left|f(z)-g(z)\right|_{z\in\gamma}&<S,
\end{align}
then we will have proved that $D\subset f(D)$ and $\deg
f=\deg g=2$ by Rouché's theorem. Indeed, given $w\in D$ $f(z)-w=0$ has
the same number of solutions as $g(z)-w=0$, which is exactly 2 counted
according with multiplicity.
\begin{figure}[!hbt]
  \begin{center}
    \psfrag{A}{\small$<S$} \psfrag{B}{\small$f(z)$}
    \psfrag{C}{\small$g(z)$} \psfrag{D}{\small$R$}
    \psfrag{E}{\small$>S$}
    \includegraphics[width=5 cm]{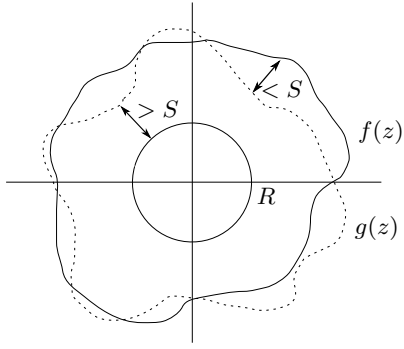}
    \caption{ Sketch of inequalities}
  \end{center}
\end{figure}
Clearly,
\begin{align}
  \nonumber \left|g(z)-w\right|_{\substack{z\in\gamma\\w\in
      D}}\geq\left|g(z)\right|_{z\in\gamma}-\left|w\right|_{w\in D}\geq
  (R^2/2-R)-R.
\end{align}
Define $S:=R^2/2-2\,R$. Since we want $S>R>0$, we require that
$R>4$. Now expand $\textrm{exp}(z/a)$ as a power series and let
$|a|=b>R$.  Then
\begin{align}
  \left|f(z)-g(z)\right|&=
  \left|\frac{z^3}{2a}+\frac{z^2}{2a}-a(z+1-a)\sum_{j=3}^\infty\frac{z^j}{j!a^j}
  \right|\leq\nonumber\\
  &\leq\frac{R^3}{2b}+\frac{R^2}{2b}+\frac{R^3}{6b^3}(3b^2e^{R/b})=
  \frac{R^2}{2b}(1+(1+e^{R/b})R).\nonumber
\end{align}
This last expression can be bounded by $\frac{R^2}{2b}(1+4R)$ as $b>R$.
Now we would like to find some $R$ such that for $b>R$, 
$\frac{R^2}{2b}(1+4R)<S$. It follows that
\[\frac{R+4R^2}{R-4}<b,\]
and this function of $R$ has a local minimum at $R\approx8.12311$. We
then conclude that given $R=8.12311$ $b$ must be larger than $65.9848$.

This way the triple $\left(f_a,D(0,R),f(D(0,R))\right)$ is
polynomial-like of degree two for $|a|\geq 66$.


\begin{remark}
Numerical experiments suggest that $|a|>10$ would be enough.
\end{remark}


\section{Semi-hyperbolic components: Proof of Theorem B}
\label{SectionSemiHyperbolicComponents}
In this section we deal with the set of parameters $a$ such that the
free singular value is attracted to a periodic orbit. We denote this
set by $H$ and it naturally splits into the pairwise disjoint subsets
\begin{align}
  H_p^v&=\{a\in\mathbb{C}\vert
  \mathcal{O}^+(\asv{a}) \textrm{ is attracted to a periodic orbit of period }p\}\nonumber \\
  H_p^c&=\{a\in\mathbb{C}\vert \mathcal{O}^+(-1) \textrm{ is attracted
    to a periodic orbit of period }p\}.\nonumber
\end{align}
where $p\geq1$. We will call these sets \emph{semi-hyperbolic
  components}.

It is immediate from the definition that semi-hyperbolic components
are open. Also connecting with the definition in the previous section
we have $H^c=\cup_{p\geq 1}H_p^c$ and $H^v=\cup_{p\geq 1}H_p^v$.

As a first observation note that, by Theorem
\ref{TheoremTheFamilyIsPolynomialLike}, every connected component of
$H_p^c$ for every $p\geq1$ is bounded. Indeed, for large values of $a$
the function $f_a(z)$ is polynomial-like and hence the critical orbit
cannot be converging to any periodic cycle, which partially proves
Theorem B, Part \ref{PartDTmaB}). We shall see that, opposite to this
fact, all components of $H_p^v$ are unbounded. We start by showing
that no semi-hyperbolic component in $H_p^c$ can surround $a=0$, by
showing the existence of continuous curves of parameter values,
leading to $a=0$, for which the critical orbit tends to
$\infty$. These curves can be observed numerically in Figure
\ref{FigureParametersEscapeToInfinityZoomCrab} in the previous
section.

\begin{proposition}
  \label{PropositionNotSurroundingZero}
  If $\gamma$ is a closed curve contained in a component $W$ of $H^c\cup
  C^c$, then $\mathrm{ind}(\gamma,0)=0$.
\end{proposition}
\begin{proof}
  We shall show that there exists a continuous curve $a(t)$ such that
  $f_{a(t)}^n(-1)\stackrel{n\to\infty}{\longrightarrow}\infty$ for all
  $t$. It then follows that $a(t)$ would intersect any curve $\gamma$
  surrounding $a=0$. But if $\gamma\subset H^c\cup C^c$, this is
  impossible. For $a\neq0$ we conjugate $f_a$ by $u=z/a$ and obtain
  the family $g_a(u)=\lambda(e^u(au+1-a)-1+a)$. Observe that
  $g_0(u)=\lambda(e^u-1)$. The idea of the proof is the following. As
  $a$ approaches 0, the dynamics of $g_a$ converge to those of $g_0$.
  In particular we find continuous invariant curves
  $\{\Gamma_k^a(t),k\in\mathbb{Z}\}_{t\in(0,\infty)}$ (Devaney hairs
  or dynamic rays) such that
  $\mathrm{Re}\,\Gamma_k^a(t)\stackrel{t\to\infty}{\longrightarrow}\infty$
  and if $z\in\Gamma_k^a(t)$ then $\mathrm{Re}\,g_a^n(z)\to\infty$.
  These invariant curves move continuously with respect to the
  parameter $a$, and they change less and less as a approaches 0,
  since $g_a$ converges uniformly to $g_0$.

  On the other hand, the critical point of $g_a$ is now located at
  $c_a=-1/a$. Hence, when $a$ runs along a half circle around 0, say
  $\eta_t=\{te^{i\alpha},\pi/2\leq\alpha\leq3\pi/2\}$, $c_a$ runs along a
  half circle with positive real part, of modulus $|c_a|=1/t$.

\begin{figure}[!hbt]
  \begin{center}
    \psfrag{A}{\small$a(t)$} \psfrag{0}{0}
    \psfrag{G0}{\small$\Gamma_0^a$} \psfrag{G1}{\small$\Gamma_1^a$}
    \psfrag{GM1}{\small$\Gamma_{-1}^a$}
    \includegraphics[height=0.40\textwidth]{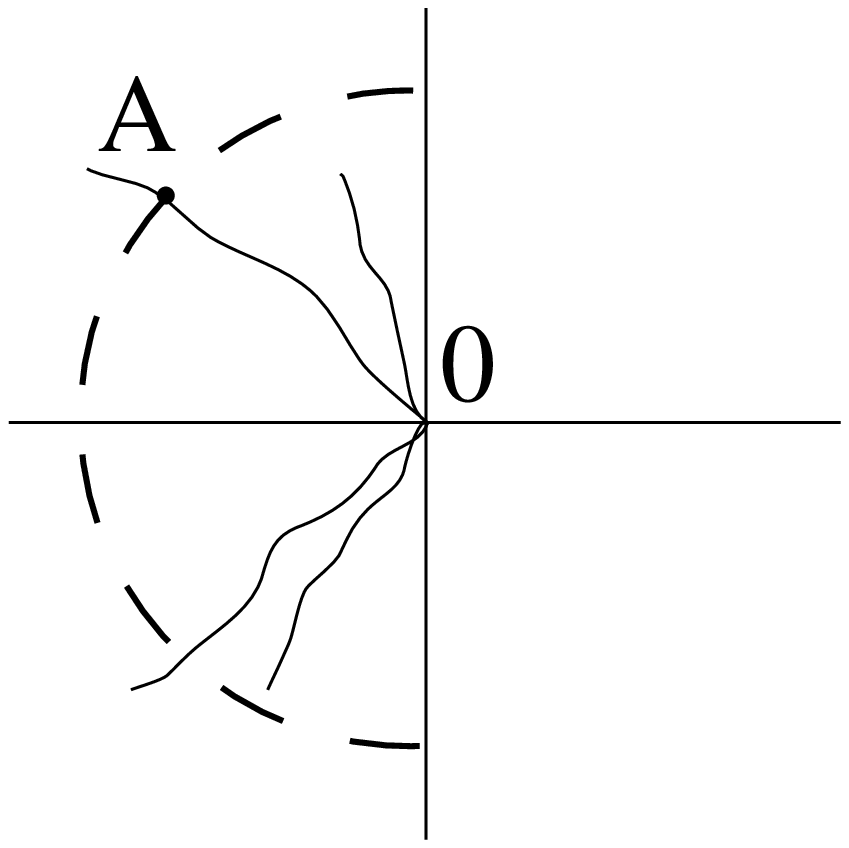}
    \includegraphics[height=0.40\textwidth]{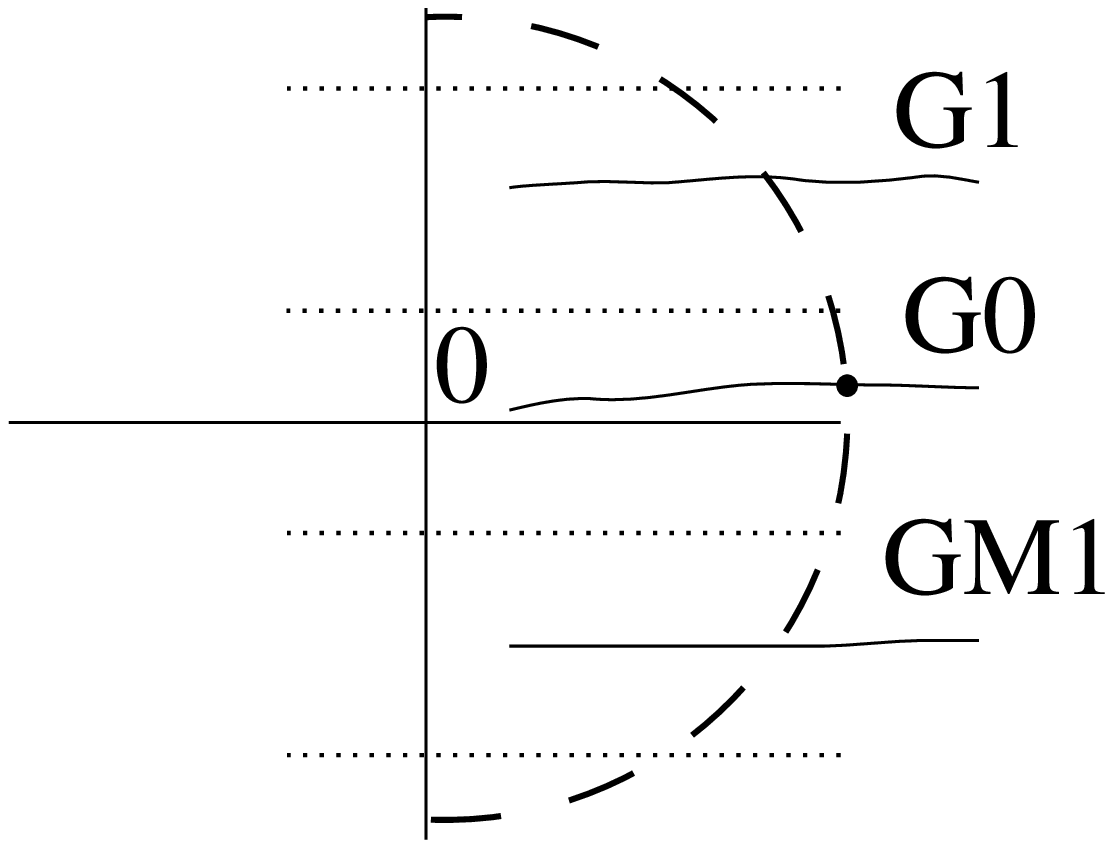}
    \caption{ {\bf Right:} Parameter plane {\bf Left: } Dynamical plane
      of $g_a(z)$.}
  \end{center}
\end{figure}

If $t$ is small enough, this circle must intersect, say, $\Gamma_0^a$
in at least one point. This means that there exists at least one
$a(t)\in\eta_t$ such that
$g_{a(t)}^n(c_{a(t)}\stackrel{n\to\infty}{\longrightarrow}\infty)$.
Using standard arguments (see for example
\cite{FagellaLimitingDynamics}) it is easy to see that we can choose
$a(t)$ in a continuous way so that $a(t)\stackrel{t\to0}{\longrightarrow}0$.
Undoing the change of variables, the conclusion follows.

\end{proof}

We would like to show now that all semi-hyperbolic components are simply connected. We first prove a preliminary lemma.


\begin{lemma}\label{LemmaHyperbolicOrbitsAreBounded}
  Let $U\subseteq H_p^v$ with $\bar U$ compact. Then there is a constant
  $C>0$ such that for all $a\in U$ the elements of the attracting
  hyperbolic orbit,  $z_j(a)$, satisfy $|z_j(a)|\leq C$, $j=1,\ldots,p$.
\end{lemma}

\begin{proof}
  If this is not the case, then for some $1\leq j\leq p,$
  $z_j(a)\to\infty$ as $a\to a_0\in \partial U$ with $a\in U$. But as
  long as $a\in U$, $z_j(a)$ is well defined, and its multiplier
  bounded (by 1). Therefore,
  \[\prod_{j=1}^p|f'_a(z_j(a))|=\prod_{j=1}^p|\lambda e^{z_j(a)/a}||z_j(a)+1|<1.\]
  Now, we claim that $z_j(a)+1$ does not converge to 0 for any $1\leq
  j\leq p$ as $a$ goes to $a_0$. Indeed, if this was the case,
  $z_j(a)$ would converge to -1, which has a dense orbit around the
  Siegel disc, but as the period of the periodic orbit is fixed, this
  contradicts the assumption. Hence $\prod_{j=1}^p|z_j(a)+1|\to\infty$
  and necessarily $\prod_{j=1}^p |e^{z_j(a)/a}|\to 0$ as $a$ goes to
  $a_0$. This implies that at least one of these elements goes to 0,
  say $|e^{z_j(a)/a}|\to0$. But this means that $z_{j+1}(a)\to\lambda
  a_0(a_0-1)=v_{a_0}$ as $a\to a_0$.  Now the first $p-1$ iterates of
  the orbit of $\asv{a_0}$ by $f_{a_0}$ are finite. Since $f_a$ is
  continuous with respect to $a$ in $\bar U$, these elements cannot be
  the limit of a periodic orbit, with one of its points going to
  infinity. In particular we would have
  $f^{p-1}_a(z_{j+1}(a))=z_j(a)\to f^{p-1}_{a_0}(\asv{a_0})$ which
  contradicts the assumption.

\end{proof}

With these preliminaries, the proof of simple connectedness is
standard (see \cite{BakerRipponExponentials} or
\cite{DevaneyExpHairs99}).

\begin{proposition}{(Theorem B, Part \ref{PartATmaB})}\label{PropositionHyperbolicSimplyConnected}
  For all $p\geq 1$ every connected component $W$ of $H_p^v$ or
  $H_p^c$ is simply connected.
\end{proposition}

\begin{proof}
  Let $\gamma\subset W$ a simple curve bounding a domain $D$. We will
  show that $D\subset W$. Let $g_n(a)= f^{np}_a(\asv{a})$ (resp. $
  f_a^{np}(-1)$).  We claim that $\{g_n\}_{n\in\mathbb{N}}$ is a
  family of entire functions for $a \in D$. Indeed, $f_a(\asv{a})$ has
  no essential singularity at $a=0$ (resp. $f_a(-1)$ has no essential
  singularity as $0\notin D$), neither do
  $f_a^n(f_a(\asv{a})),\,n\geq1$ (resp. $f_a^n(f_a(-1)),\,n\geq1$) as
  the denominator of the exponential term simplifies.

  By definition $W$ is an open set, therefore there is a neighbourhood
  $\gamma \subset U\subset W$. By Lemma
  \ref{LemmaHyperbolicOrbitsAreBounded} $|z_j(a)|<C,\,j=1,\ldots,p$
  and it follows that $\{g_n(a)\}_{n\in\mathbb{N}}$ is uniformly
  bounded in $U$, since it must converge to one point of the
  attracting cycle as $n$ goes to $\infty$. So by Montel's theorem and
  the Maximum Modulus Principle, this family is normal, and it has a
  sub-sequence convergent in $D$. If we denote by $G(a)$ the limit
  function, $G(a)$ is analytic and the mapping $H(a)=f_a^p(G(a))-G(a)$
  is also analytic. By definition of $H_p$, $H(a)$ is identically zero
  in $U$, and by analytic continuation it is also identically zero in
  $D$. Therefore $G(a)=z(a)$ is a periodic point of period $p$.

  Now let $\chi(a)$ be the multiplier of this periodic point of period
  $p$. This multiplier is an analytic function which satisfies
  $|\chi(a)|<1$ in $U$, and by the Maximum Modulus Principle the same
  holds in $D$. Hence $D\subset H_p^v$ (resp. $D\subset H^c_p$).

\end{proof}

The following lemma shows that the asymptotic value itself can not be
part of an attracting orbit.

\begin{lemma}\label{LemmaSemiHyperbolicAsymptoticComponentsHaveNoCenters}
  There are neither $a$ nor $p$ such that $f^{p}(\asv{a})=\asv{a}$ and
  the cycle is attracting.
\end{lemma}

\begin{proof}
  It cannot be a super-attracting cycle since such orbit must contain
  the critical point and its forward orbit, but the critical orbit is
  accumulating on the boundary of the Siegel disc and hence its orbit
  cannot be periodic.

  It cannot be attracting either, as the attracting basin must contain
  a singular value different from the attracting periodic point
  itself, and this could only be the critical point. But, as before,
  the critical point cannot be there. The conclusion then follows.

\end{proof}

We can now show that all components in $H_p^v$ are unbounded, which is
part of Part \ref{PartBTmaB}) of Theorem B. The proof is also
analogous to the exponential case (see \cite{BakerRipponExponentials}
or \cite{DevaneyExpHairs99}).

\begin{theorem}\label{TheoremSemiHyperbolicAsymptoticComponentsAreUnbounded}
  Every connected component $W$ of $H_p^v$ is unbounded for $p\geq1$.
\end{theorem}

\begin{proof} From Lemma \ref{LemmaHyperbolicOrbitsAreBounded} above,
  the attracting periodic orbit $z(a)$ of Proposition
  \ref{PropositionHyperbolicSimplyConnected} above is not only
  analytic in $W$ but as $\mathrm{lim\,sup}|\chi(a)|\leq1$ for $a \in
  W$, $z(a)$ has only algebraic singularities at $b\in\partial
  W$. These singularities are in fact points where $\chi(b)=1$ by the
  implicit function theorem. This entails that the boundary of $W$ is
  comprised of arcs of curves such that $|\chi(a)|=1$.

  The multiplier in $W$ is never 0 by Lemma
  \ref{LemmaSemiHyperbolicAsymptoticComponentsHaveNoCenters}, thus if
  $W$ is bounded, it is a compact simply-connected domain bounded by
  arcs $|\chi(a)|=1$. Now $\partial\chi(W)\subset\chi(\partial
  W)\subset\{\chi\vert |\chi|=1\}$ but by the minimum principle this
  implies $0\in\chi(W)$ against assumption.

\end{proof}

To end this section we show the existence of the largest
semi-hyperbolic component, the one containing a segment $[r,\infty)$
for $r$ large, which is Theorem B, Part \ref{PartCTmaB}).

\begin{theorem}\label{TheoremH1vIsUnbounded} The parameter plane of
  $f_a(z)$ has a semi-hyperbolic component $H_1^v$ of period 1 which
  is unbounded and contains an infinite segment.
\end{theorem}
\begin{proof}
  The idea of the proof is to show that for $a=r>0$ large enough there
  is a region $\mathcal{R}$ in dynamical plane such that
  $\overline{f_a(\mathcal{R})}\subset\mathcal{R}$. By Schwartz's lemma
  it follows that $\mathcal{R}$ contains an attracting fixed point. By
  Theorem \ref{TheoremTheFamilyIsPolynomialLike} the orbit of
  $\asv{a}$ must converge to it. Not to break the flow of exposition,
  the detailed estimates of this proof can be found in the Appendix.

\end{proof}

\begin{remark}
  The proof can be adapted to the case $\lambda=\pm i$ showing that
  $H_1^v$ contains an infinite segment in $i\mathbb{R}$. Observe that
  this case is not in the assumptions of this paper since $z=0$ would
  be a parabolic point.
\end{remark}




\subsection{Parametrisation of $H_p^v$: Proof of Theorem B, Part \ref{PartBTmaB}}

In this section we will parametrise connected components $W\subset
H_p^v$ by means of quasi-conformal surgery.  In particular we will
prove that the multiplier map $\chi:W\to\mathbb{D}^*$ is a universal
covering map by constructing a local inverse of $\chi$. The proof is
standard.

\begin{theorem}\label{TheoremParametrizationSemiHyperbolicAsymptotic}
  Let $W\subset H_p^v$ be a connected component of $H_p^v$ and
  $\mathbb{D}^*$ be the punctured disc. Then $\chi:W\to\mathbb{D}^*$ is the
  universal covering map.
\end{theorem}
\begin{proof}
  For simplicity we will consider $W\subset H_1^v$ in the proof.  Take
  $a_0\in W$, and observe that $f^n_a(\asv{a})$ converges to $z(a)$ as
  $n$ goes to $\infty$, where $z(a)$ is an attracting fixed point of
  multiplier $\rho_0<1$.  By Königs theorem there is a holomorphic
  change of variables \[\varphi_{a_0}:U_{a_{0}}\to\mathbb{D}\]
  conjugating $f_{a_0}(z)$ to $m_{\rho_0}(z)=\rho_{0} z$ where
  $U_{a_0}$ is a neighbourhood of $z(a_0)$.

  Now choose an open, simply connected neighbourhood $\Omega$ of
  $\rho_0$, such that $\bar\Omega\subset\mathbb{D}^*$, and for
  $\rho\in \Omega$ consider the map
  \[\xymatrix@R1.5pt{
    \psi_{\rho}:A_{\rho_0}\ar[r]^{}&A_{\rho}\\
    \hspace{2.4em}re^{i\zeta}\ar@{|->}[r]&r^\alpha
    e^{i(\zeta+\beta\log r)}, }\] where $A_r$ denotes the standard
  straight annulus $A_r=\{z\vert r<|z|< 1\}$ and


  \[\alpha=\frac{\log |\rho|}{\log |\rho_0|}, \quad
  \beta=\frac{\arg\rho-\arg\rho_0}{\log|\rho_0|}.\]


  This mapping verifies
  $\psi_\rho(m_{\rho_0}(z))=m_\rho(\psi_\rho(z))=\rho\psi_\rho(z)$. With
  this equation we can extend $\psi_\rho$ to $m_\rho(A_\rho),
  m_\rho^{2}(A_\rho),\ldots$ and then to the whole disc $\mathbb{D}$ by
  setting $\psi(0)=0$.  Therefore, the mapping $\psi_\rho$ maps the
  annuli $m_\rho^k(A_\rho)$ homeomorphically onto the annuli
  $\{z\vert |\rho^{k+1}|\leq|z|\leq \rho^k\}$.

  This mapping has bounded dilatation, as its Beltrami coefficient is
  \[\mu_{\psi_\rho}=\frac{\alpha+i\beta-1}{\alpha+i\beta+1}
  e^{2i\zeta}.\] Now define $\Psi_\rho=\psi_\rho\varphi_{a_0}$, which is
  a function conjugating $f_{a_0}$ quasi-conformally to $\rho z$ in
  $\mathbb{D}$.

  Let $\sigma_{\rho}=\Psi_\rho^*(\sigma_0)$ be the pull-back by
  $\Psi_\rho$ of the standard complex structure $\sigma_0$ in
  $\mathbb{D}$. We extend this complex structure over $U_{a_0}$ to
  $f_{a_0}^{-n}(U_{a_0})$ pulling back by $f_{a_0}$, and prolong it to
  $\mathbb{C}$ by setting the standard complex structure on those
  points whose orbit never falls in $U_{a_0}$. This complex structure
  has bounded dilatation, as it has the same dilatation as
  $\psi_\rho$. Observe that the resulting complex structure is the
  standard complex structure around 0, because no pre-image of
  $U_{a_0}$ can intersect the Siegel disc.

  Now apply the Measurable Riemann Mapping Theorem (with dependence
  upon parameters, in particular with respect to $\rho$) so we have a
  quasi-conformal integrating map $h_{\rho}$ (which is conformal where
  the structure was the standard one) so that
  $h_\rho^*\sigma_0=\sigma_\rho$. Then the mapping $g_{\rho}=h\circ
  f\circ h^{-1}$ is holomorphic as shown in the following diagram:
\begin{displaymath}
\xymatrix{
(\mathbb{C},\sigma_{\rho'}) \ar[r]^{\psi f_a\psi^{-1}} \ar[d]^{h_{\rho'}} 
& (\mathbb{C},\sigma_{\rho'}) \ar[d]^{h_{\rho'}}  \\
(\mathbb{C},\sigma_{0})\ar[r]^{g_{\rho'}}&(\mathbb{C},\sigma_{0})}
\end{displaymath}
Moreover, the map $\rho\mapsto h_\rho(z)$ is holomorphic for any given
$z\in\mathbb{C}$ since the almost complex structure $\sigma_\rho$
depends holomorphically on $\rho$. We normalise the solution given by
the Measurable Riemann Mapping Theorem requiring that -1, 0 and
$\infty$ are mapped to themselves. This guarantees that
$g_{\rho}(z)$ satisfies the following properties:
\begin{itemize}
\item $g_\rho(z)$ has 0 as a fixed point with rotation number
  $\lambda$, so it has a Siegel disc around it,
\item $g_\rho(z)$ has only one critical point, at -1 which is a simple
  critical point,
\item $g_\rho(z)$ has an essential singularity at $\infty$,
\item $g_\rho(z)$ has only one asymptotic value with one finite
  pre-image.
\end{itemize}
Moreover $g_{\rho}(z)$ has finite order by Theorem
\ref{TheoremAhlforsHolderContinuity}. Then Theorem
\ref{TheoremCharacterizationOfTheFamily} implies that
$g_{\rho}(z)=f_b(z)$ for some $b\in\mathbb{C}^*$. Now let's summarise
what we have done.

Given $\rho$ in $\Omega\subset\mathbb{D}^*$ we have a $b(\rho)\in
W\subset H_1^v$ such that $f_{b(\rho)}(z)$ has a periodic point with
multiplier $\rho$. We claim that the dependence of $b(\rho)$ with
respect to $\rho$ is holomorphic. Indeed, recall that $\asv{a}$ has
one finite pre-image, $a-1$. Hence $h_\rho(a-1)=b(\rho)-1$ which
implies a holomorphic dependence on $\rho$.

We have then constructed a holomorphic local inverse for the
multiplier. As a consequence, $\chi:H\to\mathbb{D}^*$ is a covering
map and as $W$ is simply connected by Proposition
\ref{PropositionHyperbolicSimplyConnected} and unbounded by Theorem
\ref{TheoremSemiHyperbolicAsymptoticComponentsAreUnbounded},
$\chi$ is the universal covering map.

\end{proof}




\subsection{Parametrisation of $H_p^c$: Proof of Theorem B, Part \ref{PartDTmaB}}
Let $W$ be a connected component of $H_p^c$ which is bounded and
simply connected by Theorem \ref{TheoremTheFamilyIsPolynomialLike}.
The proof of the following proposition is analogous to the case of the
quadratic family but we sketch it for completeness.

\begin{proposition}\label{PropositionParametrizationSemiHyperbolicCritical}
  The multiplier $\chi:W\to \mathbb{D}$ is a conformal isomorphism.
\end{proposition}
\begin{proof}
  Let $W^*=W\backslash\chi^{-1}(0)$. Using the same surgery
  construction of the previous section we see that there exists a
  holomorphic local inverse of $\chi$ around any point
  $\rho=\chi(z(a))\in\mathbb{D}^*$, $a\in W^*$. It then follows that
  $\chi$ is a branched covering, ramified at most over one point. This
  shows that $\chi^{-1}(0)$ consists of at most one point by Hurwitz's
  formula.

  To show that the degree of $\chi$ is exactly one, we may perform a
  different surgery construction to obtain a local inverse around
  $\rho=0$. This surgery uses an auxiliary family of Blaschke
  products. For details see \cite{BrannerFagellaQCS} or
  \cite{DouadyAsterisqueSiegelHerman}.
\end{proof}



\section{Capture components: Proof of Theorem C}
\label{SectionCaptureComponents}
A different scenario for the dynamical plane is the situation where
one of the singular orbits is eventually \emph{captured} by the Siegel
disc. The parameters for which this occurs are called capture
parameters and, as it was the case with semi-hyperbolic parameters,
they are naturally classified into two disjoint sets depending whether
it is the critical or the asymptotic orbit the one which eventually
falls in $\Delta_a$. More precisely, for each $p\geq0$ we define
\[  C=\bigcup_{p\geq0}C_p^v\cup\bigcup_{p\geq0}C_p^c,\]
where
\begin{align}
  C_p^v&=\{a\in\mathbb{C}\vert f_a^p(\asv{a})\in\Delta_a,\,p\geq0\textrm{ minimal}\nonumber\},\\
  C_p^c&=\{a\in\mathbb{C}\vert f_a^p(-1)\in\Delta_a,\,p\geq0\textrm{
    minimal}\nonumber\},
  \end{align}

Observe that the asymptotic value may belong itself to $\Delta_a$
since it has a finite pre-image, but the critical point cannot. Hence
$C_0^c$ is empty.


We now show that being a capture parameter is an open condition. The
argument is standard, but we first need to estimate the minimum size
of the Siegel disc in terms of the parameter $a$. We do so in the
following lemma.



\begin{lemma}\label{LemmaOurSiegelDiskDoNotShrinkTooMuchYoccozNotSoSmallA}
  For all $a_0\neq0$ exists a neighbourhood $V$ of $a_0$ such that
  $f_a(z)$ is univalent in $D(0,R)$.
\end{lemma}
\begin{proof}
  The existence of a Siegel disc around $z=0$ implies that there is a
  radius $R'$ such that $f_{a_0}(z)$ is univalent in $D(0,R')$. By
  continuity of the family $f_a(z)$ with respect to the parameter $a$,
  there are $R>0, \varepsilon>0$ such that $f_a(z)$ is univalent in
  $D(0,R)$ for all a in the set $\{a\vert\,|a-a_0|<\varepsilon\}$.

\end{proof}

\begin{corollary}\label{CorollaryOurSiegelDiskDoNotShrinkTooMuchYoccozNotSoSmallA}
  For all $a_0\neq0$ exists a neighbourhood $a_0\in V$ such that
  $\Delta_a$ contains a disc of radius \[\frac{C}{4R}\] where
  $C$ is a constant that only depends on $\theta$ and $R$ only depends
  on $a_0$.
\end{corollary}
\begin{proof}
  For any value of $a$ the maps $f_a(z)$ and $\tilde
  f_a(z)=\frac{1}{R}\lambda a(e^{Rz/a}(Rz+1-a)-1+a)$ are affine
  conjugate through $h(z)=R\cdot z$. For $|a-a_0|<\varepsilon$,
  $\tilde f_a(z)$ is univalent on $\mathbb{D}$, thus we can apply
  Theorem \ref{TheoremSiegelDisksDoNotShrinkYoccoz} to deduce that the
  conformal capacity $\tilde\kappa_a$ of the Siegel disc
  $\tilde\Delta_a$ is bounded from below by a constant $C=C(\theta)$.
  Undoing the change of variables we obtain
  \[R\kappa=\tilde\kappa_a\geq C(\theta)\] and therefore, by
  Koebe's 1/4 Theorem, $\Delta_a$ contains a disc of radius
  $\frac{ C(\theta)}{4R}$.
\end{proof}

\begin{theorem}[(Theorem C, Part \ref{PartATmaC})]\label{TheoremCaptureOpenCondition}
  Let $a\in C_p^v$ (resp. $a\in C_p^c$) for some $p\geq0$
  (resp. $p\geq1$) which is minimal. Then there exists $\delta>0$ such
  that $D(a,\delta)\subset C_p^v$ (resp. $C_p^c$)
\end{theorem}

\begin{proof}
  Let $b=f^p_a(\asv{a})\in\Delta_a$ (resp. $b=f^p_a(-1)\in\Delta_a$).
  Assume $b\neq0$, (the case $b=0$ is easier and will be done
  afterwards).  Define the annulus $A$ as the region comprised between
  $\overline{\mathcal{O}(b)}$ and $\partial\Delta_a$ as shown in
  Figure \ref{FigureProofOpennessCaptureDefinitionAnnulus}.

  \begin{figure}[!hbt]
    \begin{center}
      \psfrag{A}{\small $\partial\Delta_a$}
      \psfrag{B}{\small $\overline{\mathcal{O}(b)}$}
      \psfrag{C}{\small $A$}
      \includegraphics[width=4cm]{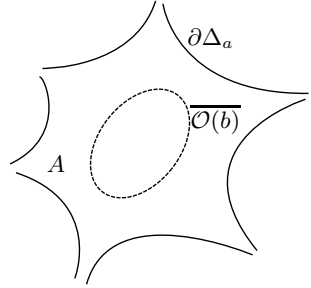}
      \caption{\label{FigureProofOpennessCaptureDefinitionAnnulus} The
        annulus $A$.}
    \end{center}
  \end{figure}

  Define $\tilde \psi$ as the restriction of the linearising
  coordinates conjugating $f_a(z)$ to the rotation
  $\mathcal{R}_\theta$ in $\Delta_a$, taking $A$ to the straight
  annulus $A(1,\varepsilon)$, where $\varepsilon$ is determined by the
  modulus of $A$. Also define a quasi-conformal mapping
  $\tilde\phi:A(1,\varepsilon )\to A(1,\varepsilon^2)$
  conjugating the rotation $\mathcal{R}_\theta$ to itself. Let $\phi$
  be the composition $\tilde \phi\circ\tilde\psi$.

  Let $\mu$ be the $f_a$ invariant Beltrami form defined as the
  pull-back $\mu=\tilde\phi^*\mu_0$ in $A$ and spread this structure
  to $\cup_nf^{-n}_a(A)$ by the dynamics of $f_a(z)$. Finally define
  $\mu=\mu_0$ in $\mathbb{C}\backslash\cup_nf^{-n}(A)$. Observe that
  $\mu=\mu_0$ in a neighbourhood of 0. Also $\phi$ has bounded
  dilatation, say $k<1$, which is also the dilatation of $\mu$.

  Now let $\mu_t=t\cdot\mu$ be a family of Beltrami forms with
  $t\in\mathbb{D}(0,1/k)$. These new Beltrami forms are integrable, since
  $\|\mu_t\|_\infty=t\|\mu\|<\frac{1}{k}k=1$. Thus by the
  Measurable Riemann Mapping Theorem we get an integrating map
  $\phi_t$ fixing 0,-1 and $\infty$, such that
  $\phi_t^*\mu_0=\mu_t$. Let $f^t=\phi_t\circ
  f_a\circ\phi_t^{-1}$,

  \[
  \xymatrix{
    (\mathbb{C},\mu_t)\ar[d]^{\phi_t}\ar[r]^{f_a}&(\mathbb{C},\mu_t)\ar[d]^{\phi_t}\\
    (\mathbb{C},\mu_0)\ar@{-->}[r]^{f^t}&(\mathbb{C},\mu_0) }
  \]

  Since $\mu_t$ is $f_a$-invariant, it follows that $f^t(z)$ preserves
  the standard complex structure and hence it is holomorphic by Weyl's
  lemma.

  Notice also that by Theorem \ref{TheoremAhlforsHolderContinuity} in
  Section \ref{SectionPreliminars} $f^t(z)$ has finite order.
  Furthermore by the properties of the integrating map and topological
  considerations, it has an essential singularity at $\infty$, a fixed
  point 0 with multiplier $\lambda$ and a simple critical point in -1.
  Finally, it has one asymptotic value $\phi_t(a)$ with one finite
  pre-image, $\phi_t(a-1)$. Hence by Theorem
  \ref{TheoremCharacterizationOfTheFamily} $f^t(z)=f_{a(t)}(z)$ for
  some $a(t)$. Now we want to prove that $a(t)$ is analytic. First
  observe that for any fixed $z\in\mathbb{C}$, the almost complex
  structure $\mu_t$ is analytic with respect to $t$. Hence, by the
  MRMT, it follows that $t\mapsto\phi_t(z)$ is analytic with respect
  to $t$. Now, $a-1$ is the finite pre-image of $\asv{a}$, so
  $\phi_t(a-1)=a(t)-1$, and this implies $a(t)=1+\phi_t(a-1)$, which
  implies that $a(t)$ is also analytic.
  
  It follows that $a(t)$ is either open or constant. But
  $f_{a(0)}=f_a$ and $f_{1}$ are different mappings since the annuli
  $\phi_0(A)=A$ and $\phi_1(A)$ have different moduli. Then $a(t)$ is
  open and therefore $\{a(t),\,t\in D(0,1/k)\}$ is an open
  neighbourhood of $a$ which belongs to $C_p^v$ (resp. $C_p^c$).
  
  If $f^p_{a_0}(\asv{a_0})=0$ (resp. $f^p_{a_0}(-1)=0$), by Lemma
  \ref{LemmaOurSiegelDiskDoNotShrinkTooMuchYoccozNotSoSmallA} and
  Corollary
  \ref{CorollaryOurSiegelDiskDoNotShrinkTooMuchYoccozNotSoSmallA}
  there exists an $\varepsilon >0$ such that for all $a$ close to
  $a_0$, $\Delta_{a_0}\supset D(0,\varepsilon)$. Hence a small
  perturbation of $f_{a_0}$ will still capture the orbit of
  $\asv{a_0}$ (resp. -1) as we wanted.

\end{proof}


\label{SectionPropertiesCaptureComponents}
The theorem above shows that capture parameters form an open set. We
call the connected components of this set, \emph{capture components},
which may be \emph{asymptotic} or \emph{critical} depending on whether
it is the asymptotic or the critical orbit which falls into
$\Delta_a$.

As in the case of semi-hyperbolic components, capture components are
simply connected. Before showing that, we also need to prove that no
critical capture component may surround $a=0$. We just state this
fact, since the proof is a reproduction of the proof of Proposition
\ref{PropositionNotSurroundingZero} above.
\begin{proposition}
  Let $\gamma$ be a closed curve in $W\subset C^v$. Then
  $\mathrm{ind}(\gamma,0)=0$.
\end{proposition}

\begin{proposition}{(Theorem C, Part
    \ref{PartBTmaC})}\label{PropositionCaptureSimplyConnected}
  All connected components $W$ of $C^v$ or $C^c$ are simply connected.
\end{proposition}

\begin{proof}
  Let $W$ be a connected component of $C^v$ or $C^c$ and
  $\gamma\subset W$ a simple closed curve. Let $D$ be the bounded
  component of $\mathbb{C}\backslash\gamma$. Let $U$ be a neighbourhood
  of $\gamma$ such that $U\subset W$. Then, for all $a\in U$,
  $f_a^n(\asv{a})$ (resp. $f_a^n(-1)$) belongs to $\Delta_a$ for
  $n\geq n_0$, and even more it remains on an invariant curve. It
  follows that $G_n^v(a)=f_a^n(\asv{a})$ (resp. $G_n^c(a)=f_a^n(-1)$)
  is  bounded in $U$ for all $n\geq n_0$.

  Since $G_n^v(a)$ is holomorphic in all of $\mathbb{C}$ (resp. in
  $\mathbb{C}^*$), we have that $G_n^v(a)$ (resp. $G_n^c(a)$) is
  holomorphic and bounded on $D$, and hence it is a normal family in
  $D$. By analytic continuation the partial limit functions must
  coincide, so there are no bifurcation parameters in $D$. Hence
  $D\subset W$.

\end{proof}

As it was the case with semi-hyperbolic components, it follows from
Theorem \ref{TheoremTheFamilyIsPolynomialLike} that all critical
capture components must be bounded, since for $|a|$ large, the
critical orbit must accumulate on $\partial \Delta_a$. This proves
Part \ref{PartCTmaC}) if Theorem C. Among all asymptotic capture
components, there is one that stands out in all computer drawings,
precisely the main component in $C_0^v$. That is, the set of
parameters for which $\asv{a}$ itself belongs to the Siegel disc.

We first observe that this component must also be bounded. Indeed, if
$\asv{a}\in\Delta_a$ then its finite pre-image $a-1$ must also be
contained in the Siegel disc. But for $|a|$ large enough, the disc is
contained in $D(0,R)$, with $R$ independent of $a$ (see Theorem
\ref{TheoremTheFamilyIsPolynomialLike}). Clearly $C_0^v$ has a unique
component, since $\asv{a}=0$ only for $a=0$ or $a=1$. This proves Part
\ref{PartDTmaC}) of Theorem C.

The ``centre'' of $C_0^v$ is $a=1$, or the map $f_a(z)=\lambda z e^z$,
for which the asymptotic value $\asv{1}=0$ is the centre of the Siegel
disc. This map is quite well-known, as it is, in many aspects, the
transcendental analogue of the quadratic family.  It is known, for
example that if $\theta$ is of constant type then $\partial \Delta_a$
is a quasi-circle and contains the critical point. This type of
properties can be extended to the whole component $C_0^v$ as shown by
the following proposition.

\begin{proposition}{(Proposition E, Part
    \ref{PartAPropositionE})}\label{PropositionHolomorphicMotionOfC0}
  If $\theta$ is of constant type then for every $a\in C_0^v$ the
  boundary of the Siegel disc is a quasi-circle that contains the
  critical point.
\end{proposition}
\begin{proof}
  For $a=1$, $f_1(z)=\lambda z e^z$ and we know that $\partial
  \Delta_a$ is a quasi-circle that contains the critical point (see
  \cite{GeyerSiegelHerman}). Define $c_n=f_1^n(-1)$, denote by
  $\mathcal{O}_a(-1)$ the orbit of -1 by $f_a(z)$ and
  \[
  \xymatrix@R1.5pt{
    H:\{c_n\}_{n\geq0}\times C_0 ^v\ar[r]&\mathbb{C}\\
    \hspace{2.2em}(c_n\phantom{,\times},\hspace{1em}a)\ar[r]
    &f_a^n(-1)}
  \]
  Then this mapping is a holomorphic motion, as it verifies
  \begin{itemize}
  \item $H(c_n,1)=c_n$, 
  \item it is injective for every $a$, as if $\asv{a}\in C_0^v$, then
    $\mathcal{O}_a(-1)$ must accumulate on $\partial \Delta_a$. Hence
    $f_a^n(-1)\neq f_a^{m}(-1)$ for all $n\neq m$.
  \item It is holomorphic with respect to $a$ for all $c_n$, an
    obvious assertion as long as $0\notin C_0^v$ which is always true.
  \end{itemize}
  Now by the second $\lambda$-lemma (Lemma \ref{LemmaSecondLambda} in
  Section \ref{SectionPreliminars}), it extends quasi-conformally to
  the closure of $\{c_n\}_{n\in\mathbb{N}}$, which contains $\partial
  \Delta_a$. It follows that for all $a\in C_0^v$, the boundary of
  $\Delta_a$ satisfies $\partial\Delta_a=H_a(\partial\Delta_a)$ with
  $H_a$ quasi-conformal, and hence $\partial\Delta_a$ is a
  quasi-circle. Since $-1\in\partial\Delta_1$, we have that
  $-1\in\partial\Delta_a$.

\end{proof}

We shall see in the next section that this same argument can be
generalised to other regions of parameter space.


 \section{Julia stability}
 \label{SectionJuliaStability}
The maps in our family are of finite type, hence $f_{a_0}(z)$ is
$\mathcal{J}$-stable if both sequences $\{f_{a}^n(-1)\}_{n\in\mathbb{Z}}$ and
$\{f_{a}^n(\asv{a})\}_{n\in\mathbb{Z}}$ are normal for $a$ in a neighbourhood of
$a_{0}$ (see \cite{McMullenComplexDynamicsAndRenormalization} or
\cite{EremenkoLyubichDynamicalEntire92}).

We define the critical and asymptotic stable components as
\begin{align}
  \mathcal {S}^c&=\{a\in\mathbb{C}\vert G_n^c(a)=f^n_a(-1)\text{ is
    normal
    in a neighbourhood of }a\},\nonumber\\
  \mathcal {S}^v&=\{a\in\mathbb{C}\vert G_n^v(a)=f^n_a(\asv{a})\text{ is
    normal in a neighbourhood of }a\},\nonumber
\end{align}
respectively. Accordingly we define critical and asymptotic unstable
components $\mathcal{U}^c$, $\mathcal{U}^v$ as their complements,
respectively. These stable components are by definition open, its
complements closed. With this notation the set of $\mathcal{J}$-stable
parameters is then $\mathcal{S}=\mathcal{S}^c\cap\mathcal{S}^v$.



Capture parameters and semi-hyperbolic parameters clearly belong to
$\mathcal{S}^c$ or $\mathcal{S}^v$. Next, we show that, because of the
persistent Siegel disc, they actually belong to both sets.
\begin{proposition}
$H^{c,v},C^{c,v}\subset\mathcal{S}$
\end{proposition}
\begin{proof}
  Suppose, say, that $a_0\in H^v$. The orbit of $\asv{a_0}$ tends to
  an attracting cycle, and hence $a_0\in\mathcal{S}^v$. In fact, since
  $H^v$ is open, we have that $a\in\mathcal{S}^v$ for all $a$ in a
  neighbourhood $U$ of $a_0$. For all these values of $a$, the critical
  orbit is forced to accumulate on $\partial\Delta_a$, hence
  $\{f_a^n(-1)\}_{n\in\mathbb{N}}$ avoids, for example, all points in
  $\Delta_a$. It follows that $\{f_a^n(-1)\}_{n\in\mathbb{N}}$ is also
  normal on $U$ and therefore $a_0\in\mathcal{S}^c$. The three
  remaining cases are analogous.

\end{proof}

Any other component of $\mathcal{S}$ not in $H$ or $C$ will be called
a \emph{queer component}, in analogy to the terminology used for the
Mandelbrot set. We denote by $Q$ the set of queer components, so that
$\mathcal{S}=H\cup C\cup Q$.

At this point we want to return to the proof of Proposition
\ref{PropositionHolomorphicMotionOfC0}, where we showed that, for
parameters inside $C_0^v$, the boundary of the Siegel disc was moving
holomorphically with the parameter. In fact, this is a general fact
for parameters in any non-queer component of the $\mathcal{J}$-stable
set.

\begin{proposition}
  Let $W$ be a non-queer component of
  $\mathcal{S}=\mathcal{S}^c\cap\mathcal{S}^v$, and $a_0\in W$. Then
  there exists a function
  $H:W\times\partial\Delta_{a_0}\to\partial\Delta_{a}$ which is a
  holomorphic motion of $\partial\Delta_{a_0}$.
\end{proposition}

\begin{proof}
  Since $W$ is not queer, we have that $W\subset H\cup C$. Let $s_a$
  denote the singular value whose orbits accumulates on
  $\partial\Delta_a$ for $a\in W$, so that $s_a\in\{-1,\asv{a}\}$. Let
  $s_a^n=f_a^n(s_a)$, and denote the orbit of $s_a$ by
  $\mathcal{O}_a(s_a)$. Then the function
  \[
  \xymatrix@R1.5pt{
    H:\mathcal{O}_{a_0}(s_{a_0})\times W\ar[r]&\mathbb{C}\\
    \hspace{2.2em}(s^n_{a_0}\phantom{,\times},\hspace{1em}a)\ar[r]
    &s_a^n}
  \]
  is a holomorphic motion, since $\mathcal{O}_a(s_a)$ must be infinite
  for all $n$, and $f_a^n(s_a)$ is holomorphic on $a$, because
  $0\notin W$. By the second $\lambda$-lemma, $H$ extends to the
  closure of $\mathcal{O}_{a_0}(s_{a_0})$ which contains $\partial
  \Delta_0$.

\end{proof}

Combined with the fact that $f_a(z)$ is a polynomial-like map of
degree 2 for $|a|>R$ (see Theorem
\ref{TheoremTheFamilyIsPolynomialLike}) we have the following
immediate corollary.

\begin{corollary}{(Proposition E, Part \ref{PartBPropositionE})}
\label{CorollaryPolynomialLikeConsequencesHolomorphicMotions}
  Let $W\subset H^v\cup C^v$ be a component intersecting $\{|z|>R\}$
  where $R$ is given by Theorem \ref{TheoremTheFamilyIsPolynomialLike}
  (in particular this is satisfied by any component of $H^v$). Then,
  \begin{enumerate}[a)]
  \item if $\theta$ is of constant type, for all $a\in W$, the
    boundary $\partial\Delta_a$ is a quasi-circle containing the
    critical point.
  \item Depending on $\theta\in\mathbb{R}\backslash\mathbb{Q}$, other
    possibilities may occur: $\partial\Delta_a$ might be a
    quasi-circle not containing the critical point, or a
    $\mathscr{C}^n,\,n\in\mathbb{N}$ Jordan curve not being a
    quasi-circle containing the critical point, or a
    $\mathscr{C}^n,\,n\in\mathbb{N}$ Jordan curve not containing the
    critical point and not being a quasi-circle. In general, any
    possibility realised by a quadratic polynomial for some rotation
    number and which persists under quasi-conformal conjugacy, is
    realised for some $f_{a}=e^{2\pi\theta i}a(e^{z/a}(z+1-a)+a-1)$.
  \end{enumerate}
\end{corollary}

\begin{remark}
  In general, for any $W\subset H^v\cup C^v$ we only need one
  parameter $a_0\in W$ for which one of such properties is satisfied,
  to have it for all $a\in W$.
\end{remark}



\appendix{Proof of Theorem \ref{TheoremH1vIsUnbounded} and numerical
  bounds}
\label{AppendixBBounds}

We may suppose $\lambda\neq\pm i$ since $\theta\neq\pm1/2$. Let
$\lambda=\lambda_1+i\lambda_2$, $\sigma=\mathrm{Sign}\left(\lambda_1\right)$ and
$\rho=\mathrm{Sign}\left(\lambda_2\right)$. We define:
\begin{figure}[hbt!]
  \begin{center}
    \psfrag{1}{$y$} \psfrag{2}{$C_2$} \psfrag{3}{$\asv{a}$}
    \psfrag{4}{$C_1$} \psfrag{5}{$s$} \psfrag{6}{$C_3$}
    \psfrag{7}{$-y$} \psfrag{8}{\tiny$f(R)$} \psfrag{9}{\large$R$}
    \includegraphics[width=6.5cm]{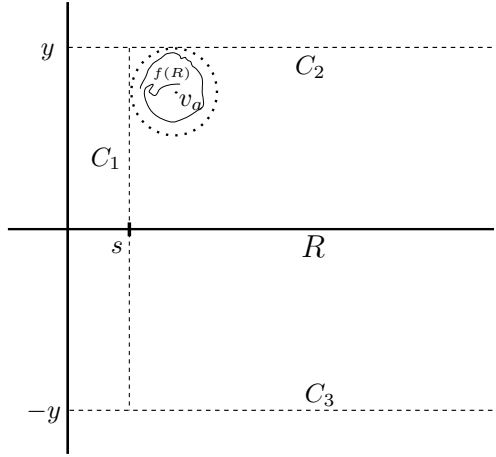}
    \caption{\label{FigureSketchCurvesInThmH1vUnbounded}Sketch of the construction
      in Thm. \ref{TheoremH1vIsUnbounded} for the case
      $\lambda_1,\lambda_2>0$.}
  \end{center}
\end{figure}
\begin{align}
  C_1:&=\{\sigma s+ti\vert |t|\leq y\}\nonumber\\
  C_2:&=\{\sigma t+i\rho y\vert t\geq s\}\nonumber\\
  C_3:&=\{\sigma t-i\rho y\vert t\geq s\}\nonumber
\end{align}
with $y>0$, $s>0$, see Figure
\ref{FigureSketchCurvesInThmH1vUnbounded} for a sketch of this
curves. Let $R$ be the region bounded by $C_1$, $C_2$, $C_3$. Recall
that $\asv{a}=\lambda(a^2-a)$ is the asymptotic value.  Note that we
will consider $a$ real, furthermore following Figure
\ref{FigureSketchCurvesInThmH1vUnbounded}, we will set $a:=-\sigma
b$ with $b>0$, as hinted by numerical experiments.  Defined this
way, the curves that are closer to $\asv{a}$ are $C_1$ and $C_2$. We
choose $y$ and $s$ in such a way that
$d(\asv{a},C_1)=d(\asv{a},C_2)$, as in Figure
\ref{FigureSketchCurvesInThmH1vUnbounded}. More precisely, \[
d(\asv{a},C_{1,2})=|\lambda_1|\left(b^2+\sigma
  b\right)-s=|\lambda_2|(b^2+\sigma b)-y
\] and hence
\[
y=\left(|\lambda_1|+|\lambda_2|\right)\left(b^2+\sigma b\right)-s.
\]
To ease notation, define $L=\left(|\lambda_1|+|\lambda_2|\right)$.
We would like some conditions over $s$ assuring that if $b>b^*$,
$d(\asv{a},f(\partial R))\leq d(\asv{a},\partial R)$, as this would
imply $f(R)\subset R$ and thus the existence of an attracting fixed
point.  We write $f_a(z)=\asv{a}+g_a(z)$ where $g_a(z)=a\cdot\lambda
e^{z/a}\cdot\left(z+1-a\right)$. Then
\[d(\asv{a}, f(\partial R))=d(0,g_a(\partial R))=|g_a(\partial
R)|.\] Therefore we need to find values such that the following
three inequalities hold
\begin{align}
  |g_a(C_1)|&< |\lambda_1|\left(b^2+\sigma b\right)-s\label{ThmH1vFirst},\\
  |g_a(C_2)|&< |\lambda_1|\left(b^2+\sigma b\right)-s\label{ThmH1vSecond},\\
  |g_a(C_3)|&< |\lambda_1|\left(b^2+\sigma
    b\right)-s\label{ThmH1vThird}.
\end{align}
For \eqref{ThmH1vFirst} to hold the following inequality needs to be
satisfied
\[b\cdot e^{-s/b}\sqrt{\left((\sigma s+\sigma b+1)+t^2\right)}
\stackrel{?}{\leq}|\lambda_1|\left(b^2+\sigma b\right)-s.\] Observe
that
\begin{align}
  b\cdot e^{-s/b}\sqrt{(\sigma s+\sigma b+1)^2+t^2}&\leq
  b\cdot e^{-s/b}\left(|\sigma(s+b)+1|+y\right)=\nonumber\\
  &=b\cdot e^{-s/b}\left(s+b+\sigma+y\right)=\nonumber\\
  &=b\cdot e^{-s/b}\left(b+\sigma+L(b^2+\sigma b)\right),\nonumber
\end{align}
so we define the following function
\[h(s)=b\cdot e^{-s/b}\left(b+\sigma+L(b^2+\sigma b)\right)-
|\lambda_1|\left(b^2+\sigma b\right)+s,\] and we will find an
argument which makes it negative. We need to find $s$ such that
$h(s)<0$ and $0<s<|\lambda_1|(b^2+\sigma b)|$.  It is easy to check
that $h(s)$ has a local minimum at
$s^*:=b\log\left(b+\sigma+L(b^2+\sigma b)\right)$ and furthermore
\[h(s^*)=b+b\log\left(b+\sigma+L(b^2+\sigma b)\right)-
|\lambda_1|\left(b^2+\sigma b\right),\] which is negative for some
$b^*$ big enough (in Appendix \ref{AppendixBBounds} we will give some
estimates on how big this $b^*$ must be as a function of $\lambda$).
This $s^*$ is again in our target interval, for a big enough $b$
(note that if $h(s^*)<0$ then $s^*<|\lambda_1|(b^2+\sigma b)|$).

From now on, let $s=s^*$, and check if \eqref{ThmH1vSecond} holds, where we
will put $s=s^*$ at the end of the calculations.
\[b\cdot e^{-\sigma t/ \sigma b}\sqrt{\left((\sigma t+\sigma
    b+1)+y^2\right)} \stackrel{?}{\leq}|\lambda_1|\left(b^2+\sigma
  b\right)-s.\] As we have done before, expand
\begin{align}b\cdot e^{-\sigma t/ \sigma b} \sqrt{\left((\sigma
      t+\sigma b+1)+y^2\right)}&\leq
  b\cdot e^{-t/b}\cdot\left(|\sigma t + \sigma b+1|+y\right)=\nonumber\\
  &=b\cdot e^{-t/b}\cdot\left( t+ b+\sigma+y\right)=\nonumber\\
  &=b\cdot e^{-t/b}\cdot\left( t+ b+\sigma +L\left(b^2+\sigma
      b\right)-s^*\right).\nonumber
\end{align}

It is easy to check that $b\cdot e^{-t/b}\cdot(b+\sigma+y)$ is a
decreasing function in $t$, and $b\cdot e^{-t/b}t$ has a local
maximum at $t=b$ and is a decreasing function for $t>b$. Then, we
can bound both terms by setting $t=s^*$, as $s^*\geq b$ whenever
$b+\sigma+L(b^2+\sigma b)$ is bigger than $e$, but this inequality
holds if all other conditions are fulfilled. Now we must only check
if
\begin{align}
  |\lambda_1|\left(b^2+\sigma b\right)-s^*&\stackrel{?}{\geq} b\cdot
  e^{-s^*/b}\cdot\left(s^*+ b+\sigma|+L
    \left(b^2+\sigma b\right)-s^*\right)=\nonumber\\
  &=b\cdot\frac{b+\sigma+L\left(b^2+\sigma b\right)}
  {b+\sigma+L\left(b^2+\sigma b\right)}=b,\nonumber
\end{align}
which is the same inequality we have for $h(s)$, thus it is also
satisfied. Inequality \eqref{ThmH1vThird} is equivalent to
\eqref{ThmH1vFirst}, hence the result follows.

Now we give numerical bounds for how big $b$ must be in Theorem
\ref{TheoremH1vIsUnbounded}. We will consider only the general case
$\lambda_1\neq0$, as the other is equivalent.

Consider the inequality
\[b\log\left(b+\sigma+L(b^2+\sigma b\right)) \leq
-b+|\lambda_1|\left(b^2+\sigma b\right)\] If this inequality holds and
$b+\sigma+L(b^2+\sigma b)>0$, we have the required estimates to
guarantee that all required inequalities in Theorem
\ref{TheoremH1vIsUnbounded} hold. The second inequality is clearly
trivial, as it holds when $b>1$.  Now, we must find a suitable $b$ for
the first.

Simplifying a $b$ factor and taking exponentials in both sides, we must check
which $b$ verify
\begin{align}
  b+\sigma+L(b^2+\sigma b) \leq
  e^{-1+|\lambda_1|\sigma}e^{|\lambda_1|b}\label{AppendixEqFindB}.
\end{align}
We can get a lower bound of $e^x$:
\begin{align}
  e^{|\lambda_1|b}\geq1+|\lambda_1|b+\frac{|\lambda_1|^2b^2}{2}
  +\frac{|\lambda_1|^3b^3}{6}.\nonumber
\end{align}
And this way if
\[b+\sigma+L(b^2+\sigma b) \leq
e^{-1+|\lambda_1|\sigma}\left(1+|\lambda_1|b+\frac{|\lambda_1|^2b^2}{2}
  +\frac{|\lambda_1|^3b^3}{6}\right),\] then is also true
\eqref{AppendixEqFindB}.  Now we must check when a degree 3 polynomial with
negative dominant term has negative values. This will be true as long
as $b>0$ is greater than the root with bigger modulus. It is
well-known (see \cite{HirstMacey97}) that a monic polynomial
$z^n+\sum_i^{n-1}a_iz^i$ has its roots in a disc of radius
$\max(1,\sum_i^{n-1}|a_i|)$, so every $b>1$ and bigger than
\begin{align}
  \frac{6}{e^{\sigma|\lambda_1|-1}|\lambda_1|^3}\cdot
  \left(
    |L-e^{\sigma|\lambda_1|-1}\frac{|\lambda_1|^2}{2}|+
    |1-e^{\sigma|\lambda_1|-1}|\lambda_1|b+L\sigma b|+
    |b+\sigma-1|
  \right)\nonumber
\end{align}
satisfies our claims.

Finer estimates for $b$ depending on $\lambda$ can
be obtained with a more careful splitting of $\lambda$ space, for instance
\begin{align}
  \{\lambda\vert\lambda\in S^1\}&=\{\lambda\in[7\pi/4,\pi/4]\}\cup
  \{\lambda\in[\pi/4,3\pi/4]\}\cup\{\lambda\in[3\pi/4,5\pi/4]\}\nonumber\\
  &\cup
  \{\lambda\in[5\pi/4,7\pi/4]\}=B_1\cup B_2\cup B_3\cup B_4.\nonumber
\end{align}
The proof can be adapted with very minor changes to this partition,
although the exposition and calculations are more cumbersome.

\section*{Acknowledgements}
The authors wish to thank Xavier Buff and Lasse Rempe for helpful
conversations. Both authors were partially supported by Spanish grant
MICINN MTM2006-05849/Consolider and European funds
MCRTN-CT-2006-035651. The second author was also partially supported
by Spanish grant MICINN MTM2008-01486 and by Catalan grant CIRIT
2005/SGR01028.

\def\cprime{$'$}

\end{document}